\documentclass{amsart}
\usepackage{amssymb, amstext, amscd, amsmath,verbatim,lineno, mathtools}
\usepackage{enumitem}
\usepackage[all,cmtip]{xy}
\usepackage{epigraph}
\usepackage{etoolbox}
\usepackage{comment}
\usepackage{latexsym}
\usepackage{amsthm}
\usepackage{geometry}
\usepackage{color}
\usepackage[all]{xy}
\usepackage{cleveref}
\usepackage[normalem]{ulem}

\newcommand{\inv}{^{-1}}
\newcommand{\z}{^{(0)}}
\newcommand{\go}{{G^{(0)}}}
\newcommand{\updag}{^\dagger}
\newcommand{\cT}{\mathcal{T}}
\newcommand{\cU}{\mathcal{U}}
\newcommand{\cF}{\mathcal{F}}
\newcommand{\cQ}{\mathcal{Q}}
\newcommand{\ZZ}{\mathbb{Z}}
\newcommand{\NN}{\mathbb{N}}
\newcommand{\TT}{\mathbb{T}}

\newcommand{\CC}{\mathbb{C}}
\newcommand{\chr}[1]{1_{#1}}
\newcommand{\chro}[1]{\tilde{1}_{#1}}
\newcommand{\chrp}[1]{1_{\psi_{#1}}}
\newcommand{\unit}[1]{{#1}^{(0)}}
\newcommand{\lineP}{\Pi}
\newcommand{\ua}{^\uparrow}

\newcommand{\ex}{\Delta}
\newcommand{\exa}{\mathcal{E}}
\newcommand{\re}{\lambda}
\newcommand{\rea}{\zeta}
\newcommand{\fe}{\re}
\newcommand{\fea}{\rea}
\newcommand{\Ooo}{O}
\newcommand{\NC}{{\operatorname{NC}}}
\newcommand{\GC}{{G_C}}
\newcommand{\GNC}{{G_{\NC}}}
\newcommand{\Stein}{{A_R(G;\Sigma)}}
\newcommand{\Steink}{{A_K(G;\Sigma)}}

\DeclareMathOperator{\ADP}{ADP}
\DeclareMathOperator{\ACP}{ACP}
\DeclareMathOperator{\AQP}{AQP}
\DeclareMathOperator{\WT}{(WT)}
\DeclareMathOperator{\i2i}{I2I}
\DeclareMathOperator{\supp}{supp}
\DeclareMathOperator{\id}{id}
\DeclareMathOperator{\spn}{span}

\DeclareMathOperator{\IM}{Im}
\DeclareMathOperator{\Char}{char}
\DeclareMathOperator{\Iso}{Iso}
\DeclareMathOperator{\interior}{Int}
\DeclareMathOperator{\bi}{BI}
\DeclareMathOperator{\su}{Su}

\DeclareMathOperator{\ord}{order}
\DeclareMathOperator{\alg}{Alg}

\newtheorem{thm}{Theorem}[section]
\newtheorem{lem}[thm]{Lemma}

\newtheorem{prop}[thm]{Proposition}

\newtheorem{corollary}[thm]{Corollary}
\newtheorem{example}[thm]{Example}
\newtheorem{defn}[thm]{Definition}
\newtheorem{ques}[thm]{Question}
\newtheorem{rmk}[thm]{Remark}

\begin{document}

\title{Intermediate Subalgebras of Cartan embeddings in rings and C*-algebras}

\author[J.H. Brown]{Jonathan H. Brown}
\address[J.H. Brown]{
Department of Mathematics\\
University of Dayton\\
300 College Park Dayton\\
OH 45469-2316 U.S.A.} \email{jonathan.henry.brown@gmail.com}

\author[L.O. Clark]{Lisa Orloff Clark}
\address[L.O. Clark]{
School of Mathematics and Statistics\\
Victoria University of Wellington\\
PO Box 600\\
Wellington 6140 N.Z.}
\email{lisa.clark@vuw.ac.nz}

\author[A.H. Fuller]{Adam H. Fuller}
\address[A.H. Fuller]{
Department of Mathematics\\
Ohio University\\
Athens\\
OH 45701 U.S.A.}
\email{fullera@ohio.edu}

\begin{abstract}
Let $D \subseteq A$ be a quasi-Cartan pair of algebras.  Then there exists a unique discrete groupoid twist $\Sigma \to G$ whose twisted Steinberg algebra is isomorphic to $A$ in a way that preserves $D$.  
In this paper, we show there is a lattice isomorphism between wide open subgroupoids of $G$ and subalgebras $C$ such that $D\subseteq C\subseteq A$ and $D \subseteq C$ is a quasi-Cartan pair. We also characterise which algebraic diagonal/algebraic Cartan/quasi-Cartan pairs have the property that every subalgebra $C$ with $D\subseteq C\subseteq A$  has $D \subseteq C$ a diagonal/Cartan/quasi-Cartan pair.   In the diagonal case, when the coefficient ring is a field, it is all of them.  Beyond that, only pairs that are close to being diagonal have this property. We then apply our techniques to C*-algebraic inclusions and give a complete characterization of which Cartan pairs $D \subseteq A$ have the property that every C*-subalgebra  $C$ with $D\subseteq C\subseteq A$ has $D \subseteq C$ a Cartan pair.   
\end{abstract}

\subjclass[2020]{ 16S99,22A22,46L05}
\keywords{Twisted groupoid, quasi-Cartan subalgebra, Cartan subalgebra}

\maketitle

\section{Introduction}

In \cite{A+8}, the authors construct a discrete groupoid twist $\Sigma \to G$ from what they call a  \emph{quasi-Cartan pair} of algebras $D\subseteq A$, see Definition~\ref{Def qC}. They show that this discrete twist is the unique twist that recovers the pair using the twisted Steinberg algebra construction introduced in \cite{A+5}, that is, $(A,D)\cong (A(G;\Sigma), A_R(\go))$.
In this way, the groupoid $G$ underlies the structure of the inclusion $D\subseteq A$.
In this paper, we study the relationship between the wide open subgroupoids of $G$ on the one hand and the subalgebras $C$ with $D\subseteq C\subseteq A$ on the other:  is there a Galois connection between them?  In Theorem~\ref{laticeIso} we indeed obtain a lattice isomorphism between wide open subgroupoids of $G$ and those intermediate  subalgebras $C$ of the inclusion $D\subseteq A$ where $D\subseteq C$ is also quasi-Cartan (Definition~\ref{Def qC}).  The question naturally arises: 
\begin{ques}
Does every intermediate subalgebra $C$ of quasi-Cartan pair $D\subseteq A$ have the property that $D\subseteq C$ is also a quasi-Cartan pair?
\end{ques}
We show in Theorem~\ref{thm: diagonal case} that this is the case if we restrict to {\em diagonal} (see Definition~\ref{Def qC}) pairs of algebras over fields.   However, the answer to the question for general quasi-Cartan pairs is dramatically no, as we explore in Section~\ref{PqC}: we show that {\em  purely quasi-Cartan inclusions}, those in which every intermediate subalgebra gives a quasi-Cartan pair, are exceedingly rare.  Corollary~\ref{cor i2i trivial} shows that if $\Sigma\to G$ is a trivial twist then $(\Stein, A_R(\go))$ is purely quasi-Cartan if and only $R$ is a field and $G$ is what we call a $\i2i$-groupoid (see Definition~\ref{def:I2I}).  As $(\Steink,A_K(\go))$ over a field $K$  is diagonal if and only if $G$ is principal and $\i2i$-groupoids are very close to being principal, purely quasi-Cartan pairs do not extend far outside the diagonal case.

Now, as most quasi-Cartan pairs $D\subseteq A$ have intermediate subalgebras $X$ with $D\subseteq X$ not quasi-Cartan, it is useful to have a condition on intermediate subalgebras which guarantees that they yield a quasi-Cartan pair.   We show in Corollary~\ref{CE qC} that if $D\subseteq C \subseteq A$ with $D\subseteq A$ Cartan and there exists a conditional expectation $\ex: A\to C$, then $D\subseteq C$ is also Cartan.



Our study is part of a larger cross-pollination between ring theory and operator algebras ($C^*$-algebras and von Neumann algebras) that has been particularly active since the introduction of Leavitt path algebras in \cite{AA05} and \cite{AMP} created a land bridge between the two subjects.  This land bridge was then expanded with the introduction of Steinberg algebras in \cite{CFST} and \cite{S}.   For example, the goal in \cite{A+8} was ``to describe it (an algebra) in terms of a concrete representation'' and includes ideas that can be traced to  theorems in von Neumann algebras \cite{FM1,FM2, V} and $C^*$-algebras \cite{Kum, Ren08}.  In fact, the broader story of reconstruction theorems like those in \cite{A+8} can be traced back and forth across the land bridge  connecting operator algebras and ring theory:  the results of \cite{Kum, Ren08} lead to the classification of shift spaces using Cuntz-Krieger algebras by Matsumoto and Matui \cite{MM}, this was extended to graph $C^*$-algebras \cite{BCW}, which was then transferred to Leavitt path algebras \cite{BCH},  then extended to graded Steinberg algebras \cite{ABHS, CR} and \cite{St2}; and then these graded Steinberg algebra results inspired the work of \cite{BFPRny} and \cite{CRST}. 

 Our present work crosses this land bridge twice.  Our initial goal was to transfer \cite[Theorem~3.3 and Theorem~4.5]{BEFPR} to the algebraic setting.  We do this in Theorem~\ref{laticeIso} and Theorem~\ref{thm: diagonal case} below.     However, we then go further and consider whether other inclusions are purely quasi-Cartan.  We show, to our surprise, that there are examples of purely quasi-Cartan inclusions that are not diagonal moreover we give a complete characterization in the non-twisted case of purely quasi-Cartan pairs in terms of the isotropy of the resulting groupoids in Corollary~\ref{cor i2i trivial}.  We then transfer these results back to the $C^*$-algebra setting in Corollary~\ref{cor pC cstar}, extending in  the work of \cite{BEFPR}.

We begin with a section of preliminaries. In Section~\ref{TSA} we provide an alternate construction of twisted Steinberg algebras using sections of the line bundle,inspired by the $C^*$-case (see, for example \cite{Ren08}) and show in Theorem~\ref{thm: sten equiv} that it is equivalent to the `usual' one from \cite{A+8}.  This construction allows us to talk about the twisted Steinberg algebra over $\Sigma\to G$ in terms of functions on $G$.  In Section~\ref{qC} we summarize the main results we use from \cite{A+8} and show how these results relate to wide subgroupoids. In Sections~\ref{SGC}--\ref{DP} we decrease our generality starting from quasi-Cartan pairs and ending with Diagonal pairs obtaining stronger results about intermediate subalgebras.   In Section~\ref{SGC}, we prove the Galois correspondence theorem showing that there is a lattice isomorphism between wide open subgroupoids and intermediate quasi-Cartan pairs.   This is a direct analog of Theorem~3.3 from \cite{BEFPR}.  In Section~\ref{CPCE} we restrict our attention to Cartan pairs and show that if $D\subseteq A$ is Cartan and $D\subseteq C\subseteq A$ and there exists a conditional expectation $\exa: A\to C$ the $D\subseteq C$ is also Cartan.  We show by example that this result does not extend to the quasi-Cartan setting. We then turn our attention to Diagonal pairs in Section~\ref{DP} and show that every intermediate subalgebra of a diagonal pair gives rise to a diagonal pair (Theorem~\ref{thm: diagonal case}) transferring  \cite[Theorem~4.5]{BEFPR} to the algebraic setting.    In Section~\ref{PqC} we prove two results Proposition~\ref{prop: two arrows}  and Proposition~\ref{bad apple} which rule out most twists from giving rise to purely quasi-Cartan pairs. We introduce $\i2i$-groupoids (Definition~\ref{def:I2I}) and show in Corollary~\ref{cor i2i trivial} that for trivial twists the pair $(\Stein, A_R(\go))$ is purely quasi-Cartan if and only if $G$ is $\i2i$ and $R$ is a field.  In Section~\ref{bimodules} we consider the more general setting of $D$-bimodules in an inclusion $A_K(\go)\subseteq \Steink$.  We show in Corollary~\ref{cor: biiff} that $A_K(\go)$-bimodules in $\Steink$ correspond to open subsets of $G$ if and only if $G$ is principal: thus attaining for $A_K(\go)$-bimodules what we expected initially for intermediate subalgebras.     Lastly, we transfer the results of Section~\ref{PqC} to the $C^*$-setting in Section~\ref{cstar}.  In an appendix, Appendix~\ref{CE}, we prove some results on conditional expectations  that are known for $C^*$-algebras but are new in the purely algebraic setting. 

In the late stages of this work we were made aware of the paper \cite{KenUrs}.
In that paper, the authors study intermediate algebras in reduced crossed-products $A \subseteq A \rtimes_r \Gamma$.
In their results, the authors use something similar to the $\i2i$ condition.
In Remark~\ref{rmk: KenUrs} we discuss the overlap of the work in \cite{KenUrs} and the results in Section~\ref{cstar}.

\subsection*{Acknowledgements} The second named author is supported by a Marsden grant from the Royal Society of New Zealand.  The initial work for this paper was done while the first named author was visiting the second named author at the Victoria University of Wellington.  The first named author would like to thank both  Prof. Clark and Prof. an Huef for their hospitality and support.

\section{Preliminaries}

Throughout this article, $R$ denotes
a commutative ring with identity, $R^\times$ is the group of units of $R$, and $A$ is an $R$-algebra.
Following \cite{A+8}, we assume that $D\subseteq A$ is a commutative subalgebra with idempotents $I(D)$.
 We further assume that $D$ is \emph{without torsion} with respect to $R$ (WT), that is
 \begin{equation*} \re e=0\implies \re=0 \quad\text{for all}\quad e\in I(D)-\{0\}, \re\in R. \end{equation*}
Condition  $\WT$ holds automatically if $R$ is a field and also if $R$ is an integral domain
and $D$ is a torsion-free R-module.  
Despite this,  $\WT$ is a restriction, as the following example shows.

\begin{example}
Suppose $R=\ZZ/(mn\ZZ)$ with $\gcd(m,n)=1$, $m,n>1$  and $D$ is a unital $R$-algebra. We claim that $D$ is not   $\WT$ with respect to $R$.    To see this, note that since $\gcd(m,n)=1$, there exists $k,\ell$ such that $mk+\ell n=-1$. Hence, $mk+1=\ell n$ for some $\ell$.
Now consider $(\ell n)^2\mod mn$.  We have
\[
(\ell n)^2=\ell n(mk+1)=\ell n \mod mn.
\]
Thus $\ell n1_D$ is an idempotent and $m( \ell n1_D)=0$, but $ m\neq 0$, proving the claim.  
\end{example}

\subsection{Groupoids}\label{groupoids}
A groupoid $G$  is a small category in which every morphism is invertible.  We identify the objects with the identity morphisms and call these the unit space, denoted $\go$.  Each element of the groupoid has a range and a source in $\go$ and we denote the maps sending an element of $G$ to its range and source by $r_G$ and $s_G$ respectively: we omit the subscript when the groupoid is clear from context.  We call a subgroupoid $H$ of $G$ {\em wide} if $\go\subseteq H$ or equivalently $H\z=\go$.  For subsets $ X$ and $Y$ of $G$ we denote
\[
XY=\{\gamma\eta: \gamma\in X, \eta\in Y, s(\gamma)=r(\eta)\}.
\]
If $X$ or $Y$ is the singleton $\{\gamma\}$ we omit the set brackets from the notation.  We say a subset $X\subseteq \go$ is {\em invariant} if $r(\gamma)\in X$ if and only if $s(\gamma)\in X$ for all $\gamma\in G$.

  We assume throughout that our groupoids are endowed with a locally compact Hausdorff topology in which multiplication and inversion are continuous.
An open bisection in a groupoid $G$ is an open set $B$ such that $r(B),s(B)$ are open in $G$ and $r|_B: B\to r(B), s|_B: B\to s(B)$ are homeomorphisms onto open subsets of $\go$.   A groupoid is {\em \'etale} if it has a basis consisting of open bisections, it is {\em ample} if it has a basis consisting of compact open bisections.
If $G$ is a Hausdorff \'etale groupoid, then the unit space $\go$ is clopen.  

We denote $\Iso(G):=\{\gamma: r(\gamma)=s(\gamma)\}$ and the interior of $\Iso(G)$ by $\interior(\Iso(G))$.  A groupoid is 
\begin{itemize}
\item {\em principal} if $\Iso(G)=\go$;
\item {\em effective} if $\interior(\Iso(G))=\go$;
\item {\em strongly effective} if every restriction of $G$ to a closed invariant subset of $\go$ is effective.
\end{itemize}

\begin{example} Let $\Gamma$ be a discrete group acting on a Cantor set $Y$ then $\Gamma\times Y$ is an ample groupoid  under the product topology where the range and source maps are given by $r(t,x)=t\cdot x$ and $s(t,x)=x$ respectively; multiplication is defined by $(t,t'\cdot x)(t',x)=(tt',x)$; and inversion is given by $(t,x)\inv=(t\inv,t\cdot x)$.
\end{example}

\begin{example} 
Let $E=(E^0,E^1,r,s)$ be a directed graph with vertices $E^0$ and edges $E^1$.  We assume $E$ is row-finite with no sources, that is  $0< |r\inv(v)|< \infty$ for all vertices $v \in E^0$.  The groupoid associated to $E$ was introduced in \cite{KPRR} and is defined as follows (using the southern convention).  A path is a sequence of edges $e_1e_2 e_3...$ such that $s(e_i)=r(e_{i+1})$.  We denote the set of finite paths by $E^*$ and the set of (one-sided) infinite paths by $E^\infty$.  For $\alpha\in E^*$ we denote the length of $\alpha$ by $|\alpha|$, that is,  the number of edges in the path $\alpha$. We say $x=x_1x_2...$ and $y=y_1y_2...\in E^\infty$ are shift equivalent with lag $k$, denoted  $x\sim_k y$, if there exists $n \in \NN$ such that $x_{n+k+i}=y_{n+i}$ for all $i\in \NN$.  Let 
\[
G_E=\{(x,k,y): x\sim_k y\} \subseteq E^{\infty} \times \mathbb{Z} \times E^{\infty}.
\]
Then $G_E$ is a groupoid with operations
\begin{align*}
& r_{G_E}(x,k,y)=x,\ s_{G_E}(x,k,y)=y,\\ 
& (x,k,y)(y,\ell,z)=(x,k+\ell, z),\ \text{and}\\
& (x,k,y)\inv=(y,-k, x).
\end{align*}
For $\alpha,\beta\in E^*$ with $s(\alpha_{|\alpha|})=s(\beta_{|\beta|})$, define the \emph{cylinder set} for $\alpha$ and $\beta$ by:
\[
Z(\alpha,\beta)=\{(\alpha x, |\alpha|-|\beta|, \beta x): x\in E^\infty, r(x)=s(\beta_{|\beta|})\},
\]
where $\alpha x$ and $\beta x$ are the infinite paths obtained by concatenation. 
The cylinder sets form a basis for a locally compact Hausdorff topology on $G_E$, with respect to which $G_E$ is an ample groupoid.
\end{example}

\subsection{Discrete twists}\label{subDT}
A discrete twist over an ample Hausdorff groupoid $G$ is a central extension of $G$ by a subgroup of units of a ring $R$ \cite[Defintion~4.1]{A+5}. We provide the details below.

\begin{defn}[{\cite[Definition~4.1]{A+5}}] \label{DT} Let  $R$ be a commutative ring with identity and let $T$ be a subgroup of $R^\times$.
We endow $R$, and hence $T$, with the discrete topology.
For an ample Hausdorff groupoid $G$, a {\em discrete twist} over $G$
is a sequence
\[
T\times \go \stackrel{i}{\to}\Sigma \stackrel{q}{\to} G
\]
where the groupoid $T\times \go$ is regarded as a group bundle with fibers $T$ and $\Sigma$ is a
groupoid with $\Sigma\z=i(\{1\}\times \go)$,  i and q are continuous groupoid homomorphisms
that restrict to homeomorphisms of unit spaces, such that the following conditions hold:
\begin{enumerate}
\item[(DT1)] The sequence is exact, in the sense that $i(T\times \{x\})=q\inv(\{x\})$ for all $x\in \go$, $i$ is injective and $q$ is surjective. 
\item[(DT2)] The groupoid $\Sigma$ is a locally trivial $G$-bundle in that for every $\gamma \in G$, there exists a compact open bisection $B \subseteq G$ containing $\gamma$ and a continuous map $\phi_{B}: B\to \Sigma$
such that 
\begin{enumerate}
\item[(i)] $q\circ \phi_{B}=\id_{B}$, and
\item[(ii)] the map $\Phi_{B}:(t,\eta)\to i(t, r(\eta))\phi_{B}(\eta)$ is a homeomorphism from $T\times B$ to $q\inv(B)$.
\end{enumerate}
\item[(DT3)] The image of $i$ is central in $\Sigma$, in the sense that $i(t,r(\sigma))\sigma=\sigma i(t,s(\sigma))$ for all
$\sigma\in \Sigma$ and $t\in T$.
\end{enumerate}
We denote a discrete twist over G by $(\Sigma,i,q)$, $\Sigma\to G$, or simply $\Sigma$ depending on the context. We
identify  $\Sigma\z$ with $\go$ via $i$ (or $q$).
\end{defn}

We remark that Definition~4.1 in \cite{A+5}, and the definition given in  \cite{A+8}, only require $G$ to be \'etale.
However, to construct Steinberg algebras,  $G$ needs to ample, so we only use discrete twists when $G$ is ample.
    Also note that, since we are viewing $R^\times$ and $T$ as discrete, we have that $\Sigma$ is ample if and only if $G$ is ample.
Moreover, while this definition is stated for a subgroup $T$ of $R^\times$, we will always take $T = R^\times$.

\section{Twisted Steinberg algebras}\label{TSA}

We develop the approach of viewing the twisted Steinberg algebra of a discrete twist $(\Sigma, i, q)$ over $G$ as sections of a line bundle, allowing us to view the elements  as functions on $G$ instead of as functions on $\Sigma$.  This is adapted for our purposes from  \cite{BFPRny}.

\subsubsection*{First description of $\Stein$: sections of a
  line bundle}
We construct $\Stein$ by considering  an $R$-bundle $L$ over $G$.  Let $L$ be the quotient
of $R\times \Sigma$ by the 
equivalence relation on $R\times \Sigma$ given by 
$(\re,\sigma)\sim (\re',\sigma')$ if and only if there
exists $t\in R^\times$ such that
$(\re',\sigma')=(t\inv \re, t\cdot \sigma)$.
Then $L$ is a Hausdorff topological space with respect to the quotient topology (of the product topology). 
Use
$[\re, \sigma]$ to denote the equivalence class of $(\re,
\sigma)$.
Observe that for any $t\in R^\times$,
\begin{equation}\label{Tact}
  [\lambda, t\cdot \sigma]=[t\lambda, \sigma].
\end{equation}
  The map $\lineP:L\rightarrow G$ is given
by
\[ \lineP:[\re, \sigma] \mapsto q(\sigma)
\]
is a continuous surjection.
There is an $R$-action on $L$ given by 
\begin{equation}
\label{eq:actionL}\re\cdot[\rea, \sigma]=[\re \rea,\sigma].\end{equation}
For  $\gamma\in G$ and
$\sigma_0\in q\inv(\gamma)$, the map $R\ni\re\mapsto
[\re, \sigma_0]\in \lineP\inv(\gamma)$ is a homeomorphism, and hence $L$ is a line bundle over $G$.  In
general, there is no canonical choice of $\sigma_0$.  However, when
$\gamma\in \unit{G}$, we have $\unit{\Sigma}\cap q^{-1}(\gamma)$ is a
singleton set,  so there is a canonical choice:
take $\sigma_0$ to be the unique element of 
$\unit{\Sigma}\cap q^{-1}(\gamma)$.  Thus, with our identification of
$\unit{\Sigma}$ and $\unit{G}$, when $x\in\unit{G}=\unit{\Sigma}$, we
sometimes identify $\lineP\inv(x)$ with $R$ via the map $\re\mapsto
[\re, x]=\re\cdot [1,x]$.

Since $\Sigma$ is locally trivial, $L$ is locally trivial as well.
Indeed, given $\ell\in L$, let $B$ be a compact open bisection of $G$
containing $\lineP(\ell)$.  Let $\phi_B:B\rightarrow \Sigma$ be a continuous map such that the map
\begin{align*}
\Phi_B \colon R^\times \times B & \longrightarrow q^{-1}(B) \\
(t,\gamma)&\longmapsto (t,r(\gamma))\cdot \phi_B(\gamma)  
\end{align*}
is a homeomorphism.
Note that such a map $\phi_B$ exists by Definition~\ref{DT}~(DT2).
For every element  $\ell_1\in
\lineP\inv(B)$, there exist unique  $\re\in R$ and $\gamma\in B$
so that  $\ell_1=[\re, \phi_B(\gamma)]$.  It follows that  the map
\begin{equation}\label{psiB}[\re,\phi_B(\gamma)]\mapsto (\re, \gamma)\end{equation} is an $R$-linear
homeomorphism $\psi_B: \lineP\inv(B) \to R\times B$.
Hence, $L$ is locally trivial.

Recall that we consider $R$ with the discrete topology (even if $R=\CC$ or another ring that typically comes with a different topology). We say a continuous function $f:G\to L$ is a \emph{section} if $\lineP(f(\gamma))=\gamma$.
 For a bisection $B \subseteq G$, we denote by $\pi_1$ the projection map $\pi_1 \colon R \times B \rightarrow R$; $(\re, \gamma) \mapsto \re$.

\begin{prop}\label{prop: locally constant}For a section $f: G\to L$, the following are equivalent:
\begin{enumerate}
\item\label{1prop} for each $\gamma\in G$ there exists a bisection $B$ containing $\gamma$ and an $R$-linear homeomorphism $\psi_B:\lineP\inv(B) \to R\times B$ such that $\pi_1\circ \psi_B \circ f|_B$ is constant;
\item\label{2prop} for every compact open bisection $B' \subseteq G$, $R$-linear homeomorphism $\psi_{B'}:\lineP\inv(B') \to R\times B'$, and $\gamma\in B'$, there exists a compact open bisection $\Ooo$ with $\gamma\in \Ooo\subseteq B'$ such that  $\pi_1\circ \psi_{B'} \circ f|_{\Ooo}$ is constant.
\end{enumerate}

\end{prop}

\begin{proof}

For \eqref{2prop} implies \eqref{1prop}, fix $\gamma \in G$ and let $B'$ be as in Definition~\ref{DT}~(DT2). By \eqref{2prop} there is an open bisection $\Ooo \subseteq B'$  such that $\pi_1\circ \psi_{B'} \circ f|_\Ooo$ is constant. Taking $B=\Ooo$ suffices. 

For \eqref{1prop} implies \eqref{2prop}, fix $B'\subseteq G$ such that there is an $R$-linear homeomorphism $\psi_{B'}:\lineP\inv(B') \to R\times B'$ and fix $\gamma \in B'$. 
 By \eqref{1prop} there is a compact open bisection $B$ in $G$ that contains $\gamma$ with $\psi_B:\lineP\inv(B) \to R\times B$ an $R$-linear homeomorphism and $\pi_1\circ \psi_B\inv \circ f\equiv \re\in R$.  We have $\{\re\}\times (B\cap B')$ is compact open in  $R\times B$.
Since $\psi_B$ and $\psi_{B'}$ are homeomorphisms, we have $\psi_{B'}\circ \psi_B\inv (\{\re\}\times (B\cap B'))$ is compact open in $R\times B'$.  Thus for each $\eta\in B'\cap B$,  there exists a compact open bisection $\Ooo_\eta$ and $\re_\eta\in R$ such that 
\[\{\re_\eta\}\times \Ooo_\eta\subseteq \psi_{B'}\circ \psi_B\inv (\{\lambda\}\times (B\cap B')).\]  
Now $ \{\{\re_\eta\}\times \Ooo_\eta: \eta\in (B\cap B')\}$  is an open cover for the compact set  $\psi_{B'}\circ \psi_B\inv (\{\re\}\times (B\cap B'))$, thus there exists a finite subcover: $ \{\{\re_{\eta_i}\}\times \Ooo_{\eta_i}\}$.   Now choose $\Ooo= \Ooo_{\eta_i}$ such that $ \gamma\in \Ooo_{\eta_i}$.
\end{proof}

We say a  section $f$  is  {\em locally constant} if $f$ satisfies the equivalent conditions of Proposition~\ref{prop: locally constant}.  That is, a section $f$ is locally constant if and only if for every compact open bisection $B$ and $R$-linear homeomorphism  $\psi_B: \lineP\inv(B) \to R\times B$, the map $\pi_1\circ \psi_B\inv \circ f$ is locally constant from $B\to R$.  For convenience, we say $f(\gamma)=0$ if   $(\pi_1\circ \psi_B\inv \circ f)(\gamma)=0$  and $f(\gamma)\neq 0$ otherwise: as $\psi_B$'s are assumed $R$-linear this does not depend on the choice of $R$-linear homeomorphism. 
Notice that a section that is locally constant is automatically continuous.
Define
\[
S_f:=\{\gamma \in G: f(\gamma)\neq 0\}.
\] 
Note that for a locally constant section $f$, $S_f$ is a clopen subset of $G$.
We say a locally constant section $f$ has \emph{compact support} if ${S_f}$ is compact in $G$.

There is a partially defined  multiplication on $L$, given by
\[
[\re, \sigma][\re',
\sigma']=[\re\re',\sigma\sigma'],\] whenever $\sigma$ and
$\sigma'$ are composable in $\Sigma$.  When $[\re,\sigma],
[\re',\sigma']\in L$ satisfy $q(\sigma)=q(\sigma')$, let
\begin{equation}\label{adds}
  [\re,\sigma]+[\re',\sigma']:=[\re +t\re',\sigma],
\end{equation}
where $t$ is the unique element of $R^\times$ so that $\sigma'=t\cdot\sigma$.  

We now let 
\[
\Stein:=\{f: G\to L: f {\text { is a compactly supported locally constant section}}\}.
\]
We endow $\Stein$ with an $R$-algebra structure where scalar multiplication is pointwise (using ~\eqref{eq:actionL}),  addition is
pointwise (using~\eqref{adds}), and multiplication is given by convolution:
\begin{equation}\label{secmul}
  f*g( \gamma)=\sum_{\eta_1\eta_2=\gamma}f(\eta_1)g(\eta_2)=\sum_{r(\eta)=r(\gamma)}
f(\eta)g(\eta\inv  \gamma).
\end{equation}

For each compact open bisection $B$ and $R$-linear homeomorphism $\psi_B:\lineP\inv(B)\to R\times B$, we can define
\[
\chrp{B}: \gamma\mapsto \begin{cases} \psi_B\inv(1,\gamma) & {\text{if }}\gamma\in B\\ 
0 &\text{otherwise}.\end{cases}
\]
Now $\Stein$ is the span of $\{\chrp{B}: B {\text{ compact open bisection}}\}$.     While $\chrp{B}$ does depend on the choice of $R$-linear homeomorphism, the next lemma shows that, locally, this choice only differs by a factor of $t$ for some $t\in R^\times$.

\begin{lem}\label{lem: characteristic} Suppose $B$ and $B'$ are compact open bisections in $G$, $\gamma\in B\cap B'$ and $\psi_B: R\times B\to \lineP\inv(B)$ and $\psi_{B'}:R\times B'\to \lineP\inv(B')$ are $R$-linear homeomorphisms.  Then there exists a compact open bisection $\Ooo_\gamma$ and a $t_\gamma\in R^\times$ such that $\chrp{B'}|_{\Ooo_\gamma}=t_\gamma\chrp{B}|_{\Ooo_\gamma}$.
\end{lem}

\begin{proof} Applying Proposition~\ref{prop: locally constant}$(2)$ to the section $\chrp{B}$ and $\gamma \in B'$, we get $\Ooo_\gamma$ such that $\pi_1\circ \psi_{B'}\inv\circ \chrp{B}|_{\Ooo_\gamma}$ is constant.
Denote this constant by $t_\gamma$.  As $\psi_{B'}(t_\gamma,\gamma)=\psi_{B}(1,\gamma)$, by the definition of the equivalence relation on $L$ we must have $t_\gamma\in R^\times$.   Now $t_\gamma\psi_{B'}(1,\eta)=\psi_{B}(1,\eta)$ for all $\eta\in \Ooo_\gamma$ and thus by definition we get $\chrp{B'}|_{\Ooo_\gamma}=t_\gamma\chrp{B}|_{\Ooo_\gamma}$ as desired.
\end{proof}

Note that if sections $f$ and $g$ are supported on bisections $B_1$ and $B_2$  respectively, and we take $\eta_i\in B_i$, then $(f*g)(\eta_1\eta_2)=f(\eta_1)g(\eta_2)$.
We can identify $A_R(\go):=\spn_R\{1_X: X \text{ compact open in } \go\}$ with a subalgebra of $\Stein$
 by
\[ A_R(\go)\to \Stein\quad\text{by}\quad \phi\mapsto
\left(\gamma\mapsto \begin{cases}
[\phi(\gamma),i(1,\gamma)] & \gamma\in \go\\ 0
&\text{otherwise}\end{cases}\right).\] Note that this identification takes
pointwise multiplication on $C_0(\go)$ to the convolution on
$C_c(G; \Sigma)$. Moreover 
\[
\phi*f(\gamma)=\phi(r(\gamma))f(\gamma)\quad\text{and}\quad f*\phi(\gamma)=f(\gamma)\phi(s(\gamma)) \quad\text{for }\quad \phi\in A_R(\go),\text{ and } f\in \Stein.
\]

\subsubsection*{Second description of $\Stein$: contravariant functions}

This description is the one found in the literature (e.g \cite{A+5}).  We say  $f:\Sigma\to R$ is locally constant if for every $\sigma\in \Sigma$, there exists a neighborhood $U$ of $\sigma$ such that $f|_U\equiv f(\sigma)$. Since we are assuming $R$ is discrete, locally constant functions are continuous.
We say $f$ is $R^\times$-contravariant, (or just contravariant when the context is clear),  if for all $t\in R^\times$ and $\sigma\in \Sigma$ we have $f(t\cdot \sigma)=t\inv f(\sigma)$.   Now define 
\[
A^c_R(G; \Sigma):=\{f: \Sigma \to R: f {\text { is a locally constant, $R^\times$-contravariant and $q(S_f)$ is compact}}\},
\]
where here $S_f\coloneqq \{\sigma \in \Sigma:f(\sigma)\neq 0\}.$
Addition and scalar multiplication in $A^c_R(G; \Sigma)$ are pointwise, and multiplication is given by a convolution formula:
\begin{equation}\label{covmul}
  f*g(\sigma)=\sum_{\eta\in G^{r(\gamma)}\cap S_f}f(\mathfrak{s}(\eta))g(\mathfrak{s}(\eta)\inv \sigma)
\end{equation} 
where $\mathfrak{s}: G\to \Sigma$ is any section of $(\Sigma,i,q)$.
It is shown in  \cite[Proposition~2.9]{A+8} that this formula does not depend on the section.

For each compact open bisection $X$ in $\Sigma$, we define $\chro{X}$ to be the contravariant function such that $q(S_{\chro{X}})= q(X)$  and $\chro{X}|_X\equiv 1$.   Then every element of $A^c_R(G;\Sigma)$ can be written as a linear combination of $\chro{X}$'s such that $q$ of the corresponding $X$'s are disjoint in $G$ by \cite[Proposition~2.8]{A+8}. 
We now show that $A^c_R(G;\Sigma)$ is isomorphic to $\Stein$. 

\begin{thm}\label{thm: sten equiv}
    Let $R$ be a commutative ring with identity and let  $(\Sigma, i, q)$ be a discrete twist over an ample Hausdorff groupoid $G$. Let $\mathfrak{s} \colon G \rightarrow \Sigma$ be any section of $(\Sigma, i, q).$
    Then the map $\Upsilon: A^c_R(G;\Sigma) \to \Stein$ defined by
    $\Upsilon(f)(\gamma)=[f(\mathfrak{s}(\gamma)),\mathfrak{s}(\gamma)]$ is an $R$-algebra isomorphism.
\end{thm}

\begin{proof}
Take any $f \in A^c_R(G;\Sigma)$.
First note that by the definition of the line bundle, $\Upsilon(f) = [f(\mathfrak{s}(\gamma)),\mathfrak{s}(\gamma)]$ does not depend on the choice of section $\mathfrak{s}$.  Further
\[
\Pi(\Upsilon(f)(\gamma)) = \Pi([f(\mathfrak{s}(\gamma)),\mathfrak{s}(\gamma)]) = q(\mathfrak{s}(\gamma)) = \gamma
\]
and hence $\Upsilon(f)$ is a section.
We have
\begin{align*}
S_{\Upsilon(f)} &= \{ \gamma \in G \colon \Upsilon(f)(\gamma) \neq 0 \} \\
&= \{ \gamma \in G \colon  [f(\mathfrak{s}(\gamma)),\mathfrak{s}(\gamma)] \neq [0, \mathfrak{s}(\gamma)] \} \\
&=\{ \gamma \in G \colon f(\mathfrak{s}(\gamma)) \neq 0 \} \\
&= q(S_f)
\end{align*}
which is compact by definition of $A_R^c(G;\Sigma)$.

We show that $\Upsilon(f)$ is locally constant. 
Fix $\gamma \in G$.
Since $f\in A^c_R(G;\Sigma)$, there exists a bisection $B$ with $\gamma \in B$, such that $f|_{\Phi_B(\{t\}\times B)}$ is constant,
where $\Phi_B \colon R^\times \times B \rightarrow q^{-1}(B)$ is the homeomorphism from Definition~\ref{DT}~(DT2).
In fact, we can assume there exists $\re\in R$ such that $f(\Phi_B((\{t\}\times B)))=t\inv \re$.
Now for $\gamma'\in B$, $\gamma'\mapsto \Phi_B(1,\gamma')$ is a section of $(\Sigma,i,q)$, so
\[
\Upsilon(f)(\gamma')=[f(\Phi_B(1,\gamma')), \Phi_B(1,\gamma')]=[ \re, \Phi_B(1,\gamma')]=[\re, (1,r(\gamma'))\phi_B(\gamma')]=[\re, \phi_B(\gamma')].
\]
Therefore, by taking $\psi_B$ as in equation~\eqref{psiB} we have $\psi_B(\Upsilon(f)(\gamma'))=(\re,\gamma')$, giving $\pi_1\left(\psi_B(\Upsilon(f)(\gamma'))\right)=\re$.
Hence, $\Upsilon(f)$ is locally constant and therefore $\Upsilon$ takes values in $\Stein$. 
By definition,  $\Upsilon$ is linear and a short computation shows it is multiplicative. 

It remains to show that $\Upsilon$ is a bijection. We do so by finding an inverse $\Theta$ of $\Upsilon$.  If $F\in \Stein$, then since $\lineP\circ F=\id_G$, for each $\sigma\in \Sigma$ we have \[\lineP([1,\sigma])\inv F(q(\sigma))=s(q(\sigma)).\]   Hence there exists $\re_\sigma \in R$ such that $\re_\sigma\cdot [1,\sigma]=F(q(\sigma))$.   
For $\sigma\in \Sigma$, define $\Theta(F)(\sigma)=\re_\sigma$.   By the definition of the line bundle we get that $\Theta(F)(t\cdot \sigma)=t\inv \Theta(F)(\sigma)$.
Also $q(S_{\Theta(F)})=S_F$, which is compact.    To see that $\Theta$ takes values in $A^c_R(G;\Sigma)$, we show that $\Theta(F)$ is locally constant. 
Fix $\sigma\in \Sigma$ and pick a compact open bisection $B$  such that there is a homeomorphism 
\[\Phi_B:  R^\times \times B\to q\inv(B)\text{ and such that }\pi_1\circ \psi_B(F)\equiv \re\]
where  $\psi_B:\lineP\inv(B)\to R\times B$ is the associated $R$-linear homeomorphism  to $\Phi_B$ as in \eqref{psiB}.   Pick $t$ such that $\sigma\in\Phi_B(\{t\}\times B)$. We show that $\Theta(F)$ is constant on $\Phi_B(\{t\}\times B)$.   Now consider $\eta\in B$ and let $\re_{t,\eta}$ be such that $[\re_{t,\eta},s(\eta)]=[1,\Phi_B(t,\eta)]\inv F(\eta)$ so that $\Theta(F)(\Phi_B(t,\eta))=\re_{t,\eta}$. 
We have 
\[F(\eta)=[\re_{t,\eta},\Phi_B(t,\eta)]=[\re_{t,\eta}, (t,r(\eta))\phi_B(\eta)]=[t\re_{t,\eta},\phi_B(\eta)].\]  So $\psi_B(F(\eta))=(t\re_{t,\eta},\eta)$ and hence $\pi_1\circ \psi_B(F(\eta))=t\re_{t,\eta}$. Thus $\re_{t,\eta}=\re t\inv$.  That is, $\Theta(F)|_{\Phi_B(\{t\}\times B)}\equiv \re t\inv$ and so $\Phi_B(\{t\}\times B)$ is a neighborhood containing $\sigma$ on which $\Theta(F)$ is constant.
That is  $F$ locally constant implies $\Theta(F)$ is locally constant. 

Now for $f\in A^c(G; \Sigma)$ we have $[1,\sigma]\inv \Upsilon(f)(q(\sigma))=[f(\sigma),s(\sigma)]$ so that $\Upsilon(f)(q(\sigma))=f(\sigma)\dot [1,\sigma]$. Thus, $\Theta(\Upsilon(f))(\sigma)=f(\sigma)$ for all $\sigma$.   Moreover $\Upsilon(\Theta(F))(\gamma)=[\re_{\mathfrak{s}(\gamma)}, \mathfrak{s}(\gamma)]=\re_{\mathfrak{s}(\gamma)}\cdot [1,\mathfrak{s}(\gamma)]=F(\gamma)$.  So we have  $\Theta=\Upsilon\inv$.  
\end{proof}

In the sequel, we will identify $\Stein$ and $A^c_R(G; \Sigma)$ and use whichever description is most convenient for our purposes.
Thus, we will drop the $c$ from the notation $A^c_R(G; \Sigma)$. 

Let $B$ be a compact open  bisection in $G$ and $\psi_B: \lineP\inv(B)\to R\times B$ an $R$-linear homeomorphism.
Let $\Theta = \Upsilon^{-1}$ as in the proof of Theorem~\ref{thm: sten equiv}.
Take $X=\psi_B\inv(\{1\}\times B)$, then 
\begin{equation}
\label{eq:xpsib}
\tilde{1}_X=\Theta(1_{\psi_B}).
\end{equation}
By  Proposition~2.8 in \cite{A+8} for $f\in \Stein$ there exists open bisections $U_1, \ldots, U_n\subseteq \Sigma$ with $\{q(U_i)\}_{i=1}^n$ mutually disjoint, such that $f=\sum_{i=1}^n r_i \chro{U_i}$ in the contravariant functions description of $\Stein$.
Transferring this result to the line bundle description as above we get the following result.

\begin{prop}[cf. {\cite[Proposition~2.8]{A+8}}]\label{prop regular} 
Let $R$ be a commutative ring with identity and let  $(\Sigma, i, q)$ be a discrete twist over an ample Hausdorff groupoid $G$.
For each $f\in \Stein$, there exists a finite collection $\cF$ of mutually disjoint compact open bisections in $G$, $R$-linear homeomorphisms $\psi_B: R\times B\to \lineP\inv(B)$, and elements $\re_B\in R$ for each $B\in \cF$ such that 
\[
f=\sum_{\cF} \re_B \chrp{B}.
\]
\end{prop}

\begin{rmk} If $\Sigma\to G$ is trivial in the sense that $\Sigma\cong R^\times \times G$ as topological groupoids, then $L\cong R\times G$ and locally constant sections of $L$ can be identified with locally constant functions $f:G\to R$.  That is, $\Stein\cong \{f: G\to R: f \text{ locally constant}\}$ and we omit $\Sigma$ from the notation and write $A_R(G)$.  

Note, it is possible for $\Sigma \cong R^\times \times G$ as {\em topological spaces} and not as {\em topological groupoids}. To see this, take $\omega: G*G \to R^\times$ such that $\omega(\gamma, \eta)\omega(\gamma\eta, \zeta)=\omega(\gamma,\eta\zeta)\omega(\eta,\zeta)$, then take $R^\times\times_\omega G$ to be $R^\times\times G$ with multiplication given by $(t,\gamma)(t',\eta)=(tt'\omega(\gamma,\eta),\gamma\eta)$.  Then $R^\times\times_\omega G$ is homeomorphic to $R^\times\times G$ but is not necessarily isomorphic to $R^\times\times G$ as a topological groupoid.
\end{rmk}

We end this section with a description of how open subgroupoids of $G$ give rise to subalgebras of $A_R(G; \Sigma)$. 
This is analogous to the C$^*$-algebra setting described in \cite[Lemma~2.7]{BFPRny}.

\begin{lem}\label{lem: subgroupoid to subalgebra}
Let $R$ be a commutative ring with identity and let $(\Sigma,i,q)$ be a discrete twist over an ample Hausdorff  groupoid  $G$.
Let $H$ be an open subgroupoid of $G$ and let $\Sigma|_H = q^{-1}(H)$.
Then, restricting the maps $i$ and $q$, we get a discrete twist 
$$ R^\times \times\unit{H} \rightarrow \Sigma|_H \rightarrow H.$$
Further $A_R(H;\Sigma|_H)$ is a subalgebra of $A_R(G;\Sigma)$.
\end{lem}

\begin{proof}
That $R^\times \times \unit{H} \rightarrow \Sigma|_H \rightarrow H$ is a twist follows from the definition.
That $A_R(H;\Sigma|_H)$ is a subalgebra of $A_R(G;\Sigma)$ follows from the fact that $H$ is an open subgroupoid of $G$.
\end{proof}

We end this section with a result that is well-known for twisted groupoid  $C^*$-algebras (see for example \cite[Theorem~5.1]{TFB3}, \cite[Lemma~9]{SW13}) and was proved  for (nontwisted) Steinberg algebras by  Clark,  Edie-Michell, an Huef, and  Sims \cite[Lemma~3.6]{CEHS}.

\begin{lem}\label{lem:exact}Let $R$ be a commutative ring with identity and let $(\Sigma,i,q)$ be a discrete twist over an ample Hausdorff groupoid  $G$. Let $X$ be a closed invariant subset of $\go$, and $\Omega=\go-X$.   Then $\Sigma X\to GX$ is a discrete twist and there exists a surjective homomorphisms of $R$-algebras $p: A_R(G;\Sigma)\to A_R(GX;\Sigma X)$ given by restriction of functions and an exact sequence

\[
A_R(G\Omega; \Sigma \Omega)\hookrightarrow A_R(G; \Sigma)\stackrel{p}{\to}A_R(GX;\Sigma X).
\]
\end{lem}

\begin{proof}  That $R^\times \times X \rightarrow \Sigma X \rightarrow GX$ is a twist follows from the definition. Since $X$ is closed $GX$ is closed in $G$, we have for every compact open bisection in $B$ in $G$ that $BX$ is compact open in $GX$.  Thus locally constant sections of $\Sigma\to G$ restrict to locally constant sections of $\Sigma X\to GX$, giving a map $p: A_R(G;\Sigma)\to A_R(GX;\Sigma X)$. A computation using the invariance of $X$ shows that $p$ is a homomorphism of $R$-algebras. The kernel of $p$ is precisely the sections supported on $G\Omega$.
\end{proof}
\section{Quasi-Cartan inclusions}\label{qC}
Throughout this section $R$ will denote a commutative ring with identity and $A$ will denote an $R$-algebra.
Further, we will always assume that $D\subseteq A$ is a commutative subalgebra of $A$ with idempotents $I(D)$, and that $D$ is without torsion with respect to  $R$.  
The {\em normalizers} of $D$ in $A$ are the elements of the set
\[
N(A,D):=\{n\in A: \exists k\in A \text{ such that } knk=k,\ nkn=n,\text{ and } nDk\cup kDn\subseteq D\}.
\]
It is shown in \cite{A+8} that $N(A,D)$ is an inverse semigroup with idempotents $I(D)$: in particular the $k$ in the definition above is unique and denoted $n\updag$.  As an inverse semigroup, $N(A,B)$ comes with a natural order structure given by 
\[ n\leq m  \Leftrightarrow \exists e\in I(D) \text{ such that } n=me.  \]
A normalizer $n$  is called a  {\em free normalizer} if $n\in D \text{ or } (n\updag n)(nn\updag)=0$.  Denote the set of free normalizers by
$N_f(A,D)$.

If $A=\Stein$ and $D=A_R(\go)$, then it is shown in the proof of Lemma~4.2 in \cite{A+8} that $\chro{X}$ is a normalizer for each compact open bisection $X$ in $\Sigma$.  Transferring this to the line bundle description we have that each $\chrp{B}$ is a normalizer for all compact open $B$ in $G$ and $R$-linear homeomorphisms $\psi_B:\lineP\inv(B)\to R\times B$.

\begin{defn}\label{def inclusions} Let $R$ be a commutative ring with identity and $A$ and $D$ be $R$-algebras with $D \subseteq A$.
\begin{enumerate}
\item We say a pair $(A,D)$ is an {\em inclusion} of $R$-algebras, or just and inclusion,  if $D$ is a subalgebra of $A$ and $I(D)$ is a set of local units for $A$.
\item  We say an inclusion $(A,D)$ is {\em{commutative}} if $D$ is a commutative $R$-algebra.
\item We say an inclusion $(A,D)$ is {\em{regular}} if $\spn I(D)=D$ and $\spn N(A,D)=A$.
\item A linear subspace $M\subseteq A$ is a {\em $D$-bimodule} if $DM\cup MD\subseteq M$.
\end{enumerate}
\end{defn}

\begin{rmk}\label{Alg is bimodule} Any intermediate subalgebra $C$ such that $D\subseteq C\subseteq A$ is  a $D$-bimodule. \end{rmk}

\begin{defn}[{\cite[Definition~3.1]{A+8}}]\label{def Conditional Expectation} Let $R$ be a commutative ring with identity, let $A$ be an $R$-algebra, and $D \subseteq A$ be a subalgebra.  A map $\ex:A\to D$ is called a {\em{conditional expectation}} if
\begin{enumerate}
\item $\ex$ is $R$-linear;
\item $\ex|_D=\id_D$; and
\item $\ex(dad')=d\ex(a)d'$ for all $a\in A, d,d'\in D$. 
\end{enumerate}
If, in addition, $(A,D)$ is a   {\em regular inclusion}, we say $\ex$ is {\em{faithful}} if the following holds:
\[\ex(na)=0 \text{ for all } n\in N(A,D) \implies a=0.\]
We say $\ex$ is {\em implemented by idempotents} if, for every $n\in N(A,D)$, there exists $e\in I(D)$ such that 
\[\ex(n)=en=ne=ene.\]
\end{defn}

\begin{defn}[{\cite[Definition~3.3]{A+8}}]\label{Def qC}
Let $R$ be a commutative ring with identity and suppose $(A,D)$ is a regular commutative inclusion of $R$-algebras. Suppose further that there exists a faithful conditional expectation $\ex:A\to D$ and that $D$ is  $\WT$ with respect to $R$.  Then we say the pair $(A,D)$ is
\begin{enumerate}
\item[($\ADP$)] an {\em algebraic diagonal pair} if $\spn N_f(A,D)=A$;
\item[($\ACP$)] an {\em algebraic Cartan pair} if $D$ is maximal abelian in $A$;
\item[($\AQP$)] an {\em algebraic quasi-Cartan pair} if $\ex$ is implemented by idempotents.
\end{enumerate}
\end{defn}
Note that $\ADP\implies \ACP\implies \AQP$ by \cite[Lemma~3.5 and 3.6]{A+8}.  
The motivating question for this paper is: given a quasi-Cartan inclusion $D \subseteq A$ and a subalgebra $C$ such that $D\subseteq C\subseteq A$, is  $D\subseteq C$ also quasi-Cartan?
Certainly $(C,D)$ is a commutative inclusion.  Moreover, if $\ex: A\to D$ is a conditional expectation implemented by idempotents, then so is $\ex|_C: C\to D$.  This is because $ N(C,D) \subseteq N(A,D)$. The difficulty in determining if $D\subseteq C$ is quasi-Cartan is showing that $(C,D)$ is regular. First, we collect some of the main results from \cite{A+8}.  

\begin{thm}[{\cite[Lemma~4.2, Proposition~7.1]{A+8}}]\label{thm: twists aqp etc}
\label{thm7.1}Let $R$ be a commutative ring with identity, let $(\Sigma, i, q)$ be a discrete twist  over an ample Hausdorff groupoid $G$, and assume $A_R(\go)$ is  $\WT$ with respect to $R$.  Then $(A_R(G;\Sigma),A_R(\go))$ is a commutative regular inclusion with a faithful conditional expectation $\ex: \Stein\to A_R(\go)$ given by restriction.  Moreover
\begin{enumerate}
\item $(\Stein,A_R(\go))$  is $\ADP$ if and only if $G$ is principal;
\item $(\Stein,A_R(\go))$  is $\ACP $ if and only if $G$ is effective;
\item $(\Stein,A_R(\go))$  is $\AQP$ if and only if $S_n$ is a bisection for all $n\in N(\Stein,A_R(\go))$ (Local Bisection Hypothesis).
\end{enumerate}
\end{thm}

The next corollary shows that Theorem~\ref{thm: twists aqp etc} passes to wide open subgroupoids.
\begin{corollary}\label{passes to open}  Let $R$ be a commutative ring with identity, let $(\Sigma, i, q)$ be a discrete twist  over an ample Hausdorff groupoid $G$ such that $A_R(\go)$ is  $\WT$ with respect to $R$.
Suppose $H$ is a wide open subgroupoid of $G$. Then we have the following:
\begin{enumerate}
\item\label{Dia} If $(A_R(G:\Sigma),A_R(\go))$ is an ADP, then so is $(A_R(H;\Sigma|_H),A_R(\go))$.
\item\label{Car} If $(\Stein,A_R(\go))$ is an ACP,  then so is $( A_R(H:\Sigma|_H),A_R(\go))$.
\item\label{qCar} If $(\Stein,A_R(\go))$ is an AQP,  then so is $(A(H;\Sigma|_H),A_R(\go)).$
\end{enumerate}
\end{corollary}

\begin{proof}
For item~\eqref{Dia}, \Cref{thm7.1} says $(\Stein, A_R(\go))$ is an ADP, then $G$ is principal.  But if $G$ is principal then so is $H$ and hence $(A_R(H;\Sigma|_H), A_R(\go))$ is an ADP by \Cref{thm7.1} again.  Similarly for item~\eqref{Car}, \Cref{thm7.1} says if $(\Stein,A_R(\go))$ is an ACP, then $G$ is effective.  Since $H$ is wide open in $G$ we  have \[\text{Int}(\text{Iso}(H))\subseteq \text{Int}(\text{Iso}(G))=\go=\unit{H}\]
and hence $H$ is effective and $(A_R(H;\Sigma|_H), A_R(\go))$ is ACP by \Cref{thm7.1} again.

 Finally, we show item~\eqref{qCar}. \Cref{thm7.1} says if $(\Stein,A_R(\go))$  is  AQP, then every normalizer in $N(A_R(\go),\Stein)$ is supported on a bisection.  It suffices to show the same is true for normalizers in $N(A_R(\go),A_R(H;\Sigma|_H))$.
 This follows from the fact that   \[N(A_R(\go),A_R(H;\Sigma|_H))\subseteq N(A_R(\go),\Stein).\qedhere \]  
\end{proof}

We now outline the construction of a discrete twist from an AQP/ACP/ADP, as described in \cite{A+8}. 
For a subset $U$ of the inverse semigroup $N(A,D)$, define the \emph{up-closure} of $U$ by 
\[
U\ua \coloneqq \{m\in N(A,D): \exists n\in U, n\leq m\}.
\]
If $U=U\ua$, we say $U$ is \emph{up-closed}.
If $U$ satisfies
\[\text{ for each } m,n \in U, \text{ there exists } r \in U \text{ such that }r\leq m, r \leq n,\] we say $U$ is \emph{down directed}.  A set that is both up-closed and down directed is called a {\em filter}.  An {\em ultrafilter} is a maximal proper filter with respect to subsets.  We denote the set of ultrafilters by $\cU(N(A,D))$, which is itself a groupoid: the inverse is defined as
\[
U^{-1} = U^\dagger = \{ u^\dagger \colon u \in U\};
\]
and two elements $U,V \in \cU(N(A,D))$ are composable if and only if $u^\dagger u v v^\dagger \neq 0$ for all $u \in U$ and $v \in V$. 
For two composable elements $U,V \in \cU(N(A,D))$, the product is given by
$$ 
UV = \{uv \colon u\in U,\ v \in V\}^\uparrow.
$$
That $UV$ is a filter is \cite[Lemma~3.5]{LawMarStei}; that $UV$ is an ultrafilter follows from \cite[Proposition~9.2.1]{LawsonBook}.

To build a twist, our $\Sigma$ will be the groupoid $\cU(N(A,D))$.
We say $U,W\in \cU(N(A,D))$ are equivalent if and only if
\[
\exists t\in R^\times\quad{\text{such that }}\quad U=tW \coloneqq \{tw \colon w \in W\}.
\]
Let $\cQ$ be the quotient map for this equivalence relation.
That is, for each $U \in \cU(N(A,D)$ the equivalence class of $U$ is 
\[\cQ(U) \coloneqq \{tU \colon t \in R^{\times}\}.\] 
Let $G=Q(\cU(N(A,D))$.
Then $G$ is a groupoid.
Both $G$ and $\Sigma$ have unit spaces isomorphic to  $\widehat{I(D)}$.
For each $n \in N(A,D)$ define the set
\[
Z(n):=\{U \in \cU(N(A,D)): n\in U\}.
\]
The sets $\{Z(n) \colon n \in N(A,D)\}$ form a basis for a topology on $\cU(N(A,D))$.
We endow $\Sigma = \cU(N(A,D))$ with this topology and $G = Q(\cU(N(A,D))$ with the quotient topology.
Then $G$ is an ample groupoid.
The following is a combination of results from \cite{A+8}.

\begin{thm}[{\cite[Theorem~6.6, Theorem~8.7]{A+8}}] \label{A+8 main} Suppose $R$ is a commutative ring with identity. Suppose further $(A,D)$ is a quasi-Cartan pair over $R$ and $D$ is  $\WT$ with respect to $R$.  Then 
\[ 
R^\times\times\widehat{ I(D)}\stackrel{\iota}{\to}\cU(N(A,D))\stackrel{\cQ}{\to} Q(\cU(N(A,D))
\] is a discrete twist and there exists an $R$-algebra isomorphism 
\[
\Lambda: A\to A_R(\cU(N(A,D));\cQ(\cU(N(A,D))
\]
such that $\Lambda$ restricts to an  $R$-algebra isomorphism $\Lambda|_D: D\to A_R(\widehat{I(D)})$.

Moreover, suppose $(\Sigma, i, q)$ is a discrete over an ample Hausdorff groupoid $G$ which satisfies the local bisection hypothesis.   Let $(\cU,\iota, \cQ)$ be the corresponding  twist of ultrafilters from the AQP $(\Stein,A_R(\Sigma\z))$. Then there exists isomorphisms of topological groupoids $\Phi$ making the following diagram commute
\begin{equation*}\xymatrix{R^\times \times\Sigma\z\ar[r]^i\ar[d]^{\Phi_{\Sigma\z} }& \Sigma
\ar[r]^q\ar[d]^{\Phi_\Sigma} & 
G\ar[d]^{\Phi_G}\\ 
R^\times \times\cU\z\ar[r]_\iota& \cU
\ar[r]_\cQ & 
\cQ(\cU).
}
\end{equation*}
\end{thm}

\begin{rmk}\label{Phi computation} The map $\Phi_\Sigma$ is defined by $\Phi_\Sigma(\sigma)=\left (\chro{X}U_{s(\sigma)}\right)\ua$ where $X$ is any compact open bisection in $\Sigma$ containing $\sigma$ and $U_{s(\sigma)}=\{e\in I(B): e(s(\sigma))\neq 0\}$ is the principal ultrafilter in $I(B)$ associate to $s(\sigma)$.
It is  shown in \cite[Proposition~8.4]{A+8} that this does not depend on the choice of $X$.

Note that \cite{A+8} uses the contravariant functions view.  To translate this to the line bundle view, choose a compact open  bisection $B$ in $G$ containing $q(\sigma)$ and a $R$-linear homeomorphism $\psi_B: \lineP\inv(B)\to R\times B$ such that $\psi_B(\sigma)=(1,q(\sigma))$.  Then $\Phi_\Sigma(\sigma)=\left (\chrp{B}U_{s(\sigma)}\right)\ua$.  We claim that 
\begin{equation}
\label{eqChi}\Phi_G(\gamma)=\{n\in N(\Stein, A_R(\go)): n(\gamma)\neq 0\}.\end{equation}
To see this we need to show first that for  $n\in N(\Stein, A_R(\go))$ with $n(\gamma)\neq 0$ we have $n\in \left (\chrp{B}U_{s(\sigma)}\right)\ua$ for some compact open bisection $B \subseteq G$.  But this follows from equation~\eqref{eq:xpsib} and
\cite[Proposition~8.6]{A+8}.  It remains to show that if $\chrp{B}(\gamma)$ and $\chrp{B'}(\gamma)$ are both nonzero then there exists $t\in R^\times$ such that 
\[\left (\chrp{B'}U_{s(\sigma)}\right)\ua= t\left (\chrp{B}U_{s(\sigma)}\right)\ua.\]
Since $\chrp{B}(\gamma), \chrp{B'}(\gamma)\neq 0$ we get $\gamma\in B\cap B'$.  
By Proposition~\ref{prop: locally constant}, there exists $\Ooo$ with $\gamma\in \Ooo\subseteq B\cap B'$ such that $\pi_1\circ \psi_{B'}\circ \chrp{B}|_\Ooo$ is constant $t$. Then $\chrp{B'}|_\Ooo=t\chrp{B}|_\Ooo=t\chr{\psi_{B}|_\Ooo}$.  Since $\chr{\psi_{B'}|_\Ooo}$ is a normalizer we must have $t\in R^\times$.  Lastly, invoking \cite[Proposition~8.4]{A+8} again, we get that 
\[
\left (\chrp{B'}U_{s(\sigma)}\right)\ua=\left (\chr{\psi_{B'}|_\Ooo}U_{s(\sigma)}\right)\ua=\left (t\chr{\psi_{B}|_\Ooo}U_{s(\sigma)}\right)\ua=t\left (\chrp{B}U_{s(\sigma)}\right)\ua
\]
as desired.  As each $\left (\chrp{B}U_{s(\sigma)}\right)\ua$ only contains normalizers that are nonzero at $\gamma$, we get that  $Q(\left (\chrp{B}U_{s(\sigma)}\right)\ua)=\{n\in N(A,B): n(\gamma)\neq 0\}$ and equation~\eqref{eqChi} holds.
\end{rmk}

\section{Galois correspondence}\label{SGC}

Suppose $R$ is a commutative ring with identity and $(A,D)$ is a quasi-Cartan pair.  We can thus assume that $A=\Stein$ and $D=A_R(\go)$ by 
\cite[Theorem~6.6]{A+8}.  In this section we give a Galois correspondence between
 quasi-Cartan subalgebras $C$ such that
 \[D=A_R(\go)\subseteq C\subseteq \Stein=A \]
 and wide open subgroupoids of $G$.
 We first study the lattice structure of the subalgebras $C$. For a set $\Omega \subseteq A$, we denote the algebra generated by $\Omega$, by $\alg \Omega$.

 \begin{lem}\label{intermediate lattice}
 Let $R$ be a commutative ring with identity. Suppose $D \subseteq A$ is an algebraic quasi-Cartan inclusion over a ring $R$.
 Let $\mathfrak{L}$ be the collection subalgebras $C$ with $D \subseteq C \subseteq A$ such that $D \subseteq C$ is quasi-Cartan.
 Then $\mathfrak{L}$ is a lattice with meet and join operations given by
 \begin{align*}
    C_1 \wedge C_2 &= \spn N(C_1\cap C_2, D) \text{ and }\\
    C_1 \vee C_2 &= \alg\{C_1\cup C_2\}.
 \end{align*}
 for each $C_1,\ C_2 \in \mathfrak{L}$.
 \end{lem}

 \begin{proof}
     Let $\ex:A \to D$ be the faithful conditional expectation implemented by idempotents.
     For any subalgebra $C$ such that $D \subseteq C \subseteq A$, we have that $\ex|_C$ is faithful and implemented by idempotents. 
     Thus, to show that a subalgebra $C$ such that $D \subseteq C \subseteq A$ is in $\mathfrak{L}$, it suffices to show that the inclusion $D \subseteq C$ is regular.
     
     Take $C_1, C_2 \in \mathfrak{L}$.
     Then $D \subseteq \spn N(C_1\cap C_2, D)$ is a regular inclusion by definition.
     Thus $C_1 \wedge C_2 \in \mathfrak{L}$ and $C_1 \wedge C_2 \subseteq C_1, C_2$.
     Suppose that $C\in \mathfrak{L}$ satisfies $C \subseteq C_1, C_2$.
     Then $N(C,D) \subseteq N(C_1, D) $ and $N(C,D) \subseteq N(C_2, D)$. Thus $N(C,D) \subseteq N(C_1 \cap C_2, D)$ and hence $C \subseteq C_1 \cap C_2$.  It follows that $C_1\wedge C_2$ greatest lower-bound for $C_1$ and $C_2$ in $\mathfrak{L}$.
     
    Now let $E$ be the algebra generated by $C_1 \cup C_2$. Then $E$ is the smallest algebra containing both $C_1$ and $C_2$.  It remains to show $E \in \mathfrak{L}$. 
    As $C_1 = \spn N(C_1,D)$ and $C_2 = \spn N(C_2,D)$, it follows that $E$ is spanned by elements of the form $n_1 n_2 \ldots n_k$ where each $n_i \in N(C_1,D) \cup N(C_2,D)$. As $N(C_1,D) \cup N(C_2,D) \subseteq N(A,D)$ and $N(A,D)$ is a semigroup, we have that $E = \spn N(E,D)$ and hence $E \in \mathfrak{L}.$
    \end{proof}

Given an AQP $(A,D)$, the groupoid $\Sigma = \cU(N(A,D)$ from Theorem~\ref{A+8 main} (\cite[Theorem~6.6]{A+8}) is a groupoid of ultrafilters associated to the inverse semigroup of normalizers.  We establish some general results about ultrafilters for use in the sequel. 

\begin{lem}\label{ultraSubsemigroups}
Let $T$ be an inverse semigroup with inverse $\dag$, and let $S \subseteq T$ be a  a sub-inverse semigroup.
Suppose $t^\dagger  t\in S$ for all $t\in T$.    Then
\begin{enumerate}
\item\label{idempotents} the idempotents coincide, that is, $I(S)=I(T)$;
\item\label{uparrow} if $W$ a filter in $S$, then $W^\uparrow\cap S=W$ and $W^\uparrow$ is a filter in $T$;
\item\label{TtoSfilter} if $U_T$ is a filter in $T$ and $U_T\cap S\neq \emptyset$, then $U_T\cap S$ is a filter in $S$ and $(U_T\cap S)^\uparrow =U_T$;
\item\label{TtoS} if $U_T$ is an ultrafilter in $T$ and $U_T\cap S\neq \emptyset$, then $U_T\cap S$ is an ultrafilter in $S$ and $(U_T\cap S)^\uparrow =U_T$;
\item\label{StoT} if $U_S$ is an ultrafilter in $S$, then $U_S^\uparrow$ is an ultrafilter in $T$.
\end{enumerate}
\end{lem}

\begin{proof} 
For \eqref{idempotents}, note that $I(T)=\{t^\dagger  t: t\in T\}$.  Thus $I(T)\subseteq S$ by assumption and since these are idempotents $I(T)\subseteq I(S)$.  The reverse inclusion is trivial.

For \eqref{uparrow}, let $W$ be a filter in $S$.  Since $W\subseteq W^\uparrow$, we have $W=W\cap S\subseteq W^\uparrow\cap S$.  Now if $x\in W^\uparrow\cap S$, then there exists $y\in W$ such that $y\leq x$.  But since $W$ is a filter in $S$ we get $x\in W$.
That $W^\uparrow$ is a filter is straightforward.

For \eqref{TtoSfilter}, suppose $U_T$ is a filter in $T$ and $U_T\cap S\neq \emptyset$.  We first show that $U_T\cap S$ is a filter in $S$.  To show that $U_T\cap S$ is down directed, suppose $x,y\in U_T\cap S$.  Since $U_T$ is a filter there exists $z\in U_T$  such that $z\leq x,y$.  Thus there exists an idempotent $e\in I(T)=I(S)$ (by \eqref{idempotents}) such that $z=ex$.  Since $e,x\in S$ we have $ z\in S$ so $z\in U_T\cap S$.  To see that $U_T\cap S$ is up-closed suppose $s\in U_T\cap S$ and $x\in S$ such that $s\leq x$. Since $U_T$ is a filter, $x\in U_T$ so $x\in U_T\cap S$.  

Next we show that $(U_T\cap S)^\uparrow =U_T$.  We have $(U_T\cap S)^\uparrow\subseteq U_T$ since $U_T$ is up-closed.   Suppose $x\in U_T$ and pick $y\in U_T\cap S$.  Since $U_T$ is down directed there exists $z\in U_T$ with $z\leq x, y$.  Thus, there exists an idempotent $e\in I(T)=I(S)$ (by \eqref{idempotents} again) such that $z=ex\in S$.  Thus by definition $x\in (U_T\cap S)^\uparrow$.  

For \eqref{TtoS}, suppose $U_T$ is an ultrafilter.   Fix a proper filter $W$ in $S$ such that  $U_T\cap S\subseteq W$.  Then applying \eqref{TtoSfilter} we have
\[U_T=(U_T\cap S)^\uparrow\subseteq W^\uparrow.\]   
Since $U_T$ is an ultrafilter,  $U_T=W^\uparrow$.  Now by \eqref{uparrow}, 
$W=W^\uparrow\cap S=U_T\cap S$ and hence $U_T\cap S$ is an ultrafilter. 

For \eqref{StoT},  that $U_S^\uparrow$ is a filter follows from \eqref{uparrow}.
To see it is an ultra filter, suppose $W_T$ is a proper filter in $T$ such that 
\[
U_S^\uparrow\subseteq W_T.
\]
Then $U_S=U_S^\uparrow\cap S\subseteq W_T\cap S$.  Since $U_S$ is an ultrafilter and $W_T \cap S$ is a filter by \eqref{TtoSfilter}, we have $U_S=W_T \cap S$.  
Now 
\[U_S^\uparrow = (W_T \cap S)^\uparrow =W_T\]
by \eqref{TtoSfilter}.
\end{proof}

\begin{lem}
\label{lem:upcloseup}
 Let $S \subseteq T$ be such that $S$ is an  sub-inverse semigroup of an inverse semigroup $T$ with inverse $\dagger$.  Suppose $t^\dagger  t\in S$ for all $t\in T$, and $0 \in T$.
 If $U$ and $V$ are ultrafilters in $S$ satisfying 
 $$ u^\dagger uvv^\dagger \neq 0$$
 for all $u \in U$ and $v \in V$, then
 \[ (UV)^\uparrow  = (U^\uparrow V^\uparrow)^\uparrow.\]
\end{lem}

\begin{proof}
By \cite[Lemma~3.5]{LawMarStei} and \cite[Proposition~9.2.1]{LawsonBook}, $(UV)^\uparrow \cap S$ is an ultrafilter in $S$.
Hence, $(UV)^\uparrow$ is an ultrafilter in $T$ by Lemma~\ref{ultraSubsemigroups}~(5).
Take $w \in (UV)^\uparrow$. Then, there exists $u \in U$ and $v \in V$ such that $uv \leq w$.
Note that $u \in U^\uparrow$ and $v \in V^\uparrow$, and so $w \in (U^\uparrow V^\uparrow)^\uparrow$.
Thus $(UV)^\uparrow \subseteq (U^\uparrow V^\uparrow)^\uparrow$.
As $(UV)^\uparrow$ is an ultrafilter, we have $W = (U^\uparrow V^\uparrow)^\uparrow$.
\end{proof}

In what follows, we use the following notation: 
\begin{align}
\NC&=\spn N(C,D);\\
\GC&=\{\gamma\in G: \exists c\in C, c(\gamma)\neq 0\};\\
\GNC&=\{\gamma\in G:\exists n\in N(C,D), n(\gamma)\neq 0\}=\{\gamma\in G: \exists c\in\NC, c(\gamma)\neq 0\}.
\end{align}

The following is a consequence of the previous lemmas.  If $S$ is a sub-inverse semigroup, we write $\cU(S)$ for the set of ultrafilters in $S$.

\begin{lem}\label{groupoid inclusions}  Let $R$ be a commutative ring with identity and let $(\Sigma,i,q)$ be a discrete twist over an ample Hausdorff groupoid $G$.  Suppose for  $D=A_R(G^{(0)})$ and $A=\Stein$ that $(A,D)$ is an AQP, in particular  $D$ is  $\WT$ with respect to $R$.
Let $C$ be a subalgebra such that 
$D\subseteq C\subseteq A.$ 
Denote by $\cT$ the image of $U\mapsto U^\uparrow$ from $\cU(N(C,D))$ to $\cU(N(A,D))$. Then $\cT$ is a wide open subgroupoid of $\Sigma$, 
\[\GNC=\{\gamma\in G: \exists n\in N(C,D), n(\gamma)\neq 0\}\]
is a  wide open subgroupoid of $G$ such that $\cT \cong \Sigma|_{\GNC} (=q^{-1}(\GNC))$, and \[\NC= A_R(\GNC;\Sigma|_{\GNC}).\]
\end{lem}

\begin{proof}
Since $D\subseteq A$ is quasi-Cartan, $\Sigma\cong \cU(N(A,D))$ and $G\cong Q(\cU(N(A,D)))$ by \cite[Corollary~8.10]{A+8}.  
Consider the inclusion 
\[D\subseteq \NC=\spn N(C,D).\]
Since $D\subseteq A$ is quasi-Cartan, we have $D\subseteq \NC$ is also quasi-Cartan.  Thus by Theorem~\ref{A+8 main} (\cite[Theorem 6.6]{A+8}), we have $\NC\cong A_R(\cU(N(C,D)),Q(\cU(N(C,D)))$.   

Notice that $N(C,D)$ is a sub-inverse semigroup of $N(A,D)$ such that 
for all $n \in N(A,D)$, we have 
\[n^\dagger n \in I(D) \subseteq  N(C,D).\]
Thus Lemma~\ref{ultraSubsemigroups}~\eqref{StoT} tells us that for each $U \in \cU(N(C,D))$, $U^\uparrow$ is in $\cU(N(A,D)).$
Hence $U\mapsto U^\uparrow$ maps $\cU(N(C,D))\to \cU(N(A,D))$.  Let $H\coloneqq Q(\cT) \subseteq G$.
Fix $W \in \cU(N(C,D))$ and $n \in W$.  Then 
$n\in N(C,D)\subseteq N(A,D)$.    
By applying Lemma~\ref{ultraSubsemigroups} \eqref{TtoS} we see that the open set
\[Z(n)=\{U \in \cU(N(A,D)) \colon  n\in U\} =
\{U^\uparrow\colon U \in \cU(N(C,D)), n \in U\},\]
is inside $\cT$. Thus $W^\uparrow \in Z(n) \subseteq \cT$ and hence 
$\cT$ is open in $\Sigma$.  Since $\cQ$ is an open map, 
$H$ is also open in $G$. 

Note that, for any $U \in \cU(N(C,D))$, 
$(U^\uparrow)^{-1} = (U^{-1})^\uparrow$,
and hence $\cT$ is closed under inverses.
If $U,V \in \cT$ are composable elements of $N(A,D)$ then $UV \in \cT$ by Lemma~\ref{lem:upcloseup}.
So $\cT$ is an open subgroupoid of $\Sigma$ and hence $H$ is also an open subgroupoid of $G$. 
That $\cT$ is wide follows from $D$ being a subset of $C$ and
\[\unit{\Sigma}= \{U \in \Sigma: U \cap I(D) \neq \emptyset\}.\]

Next we show $H \cong \GNC$.  Recall from Remark~\ref{Phi computation} that the isomorphism $\Phi$ from $G$ to $Q(\cU(A,D))$ 
is given by $\Phi_G(\gamma)=\{n\in N(A,D): n(\gamma)\neq 0\}$.  So 
\begin{align*}
\gamma\in H &\Leftrightarrow \text{ there exists } n\in N(C,D) \text{ with } n\in \Phi_G(\gamma) \\
&\Leftrightarrow \text{ there exists } n\in N(C,D) \text{ with } n(\gamma)\neq 0,
\end{align*}
giving $H \cong\GNC$.  Thus $\cT \cong \Sigma|_\GNC$ and we have $A_R(\Sigma|_\GNC;\GNC)\cong A_R(\cU(N(A,C))$, and $\cU(N(A,C)))\cong\NC=\spn N(A,C)$. 
\end{proof}

\begin{thm}\label{laticeIso} Let $R$ be a commutative ring with identity and $(\Sigma,i, q)$ be a discrete twist over an ample Hausdorff groupoid $G$. Suppose that $A_R(\go)\subseteq \Stein$ is quasi-Cartan.  Then \[C\mapsto G_C:=\{\gamma: \exists c\in C, c(\gamma)\neq 0\}\] defines a lattice isomorphism between subalgebras $C$ such that $A_R(\go)\subseteq C\subseteq A_R(G)$ with  $A_R(\go)\subseteq C$ quasi-Cartan and with the lattice of wide open subgroupoids of $G$ under inclusion.
\end{thm}

\begin{proof}
Suppose $A_R(\go)\subseteq C\subseteq \Stein$ with  $A_R(\go)\subseteq C$ quasi-Cartan.  Then $C=\spn N(C,D)$ by assumption and  then Lemma~\ref{groupoid inclusions} gives that $\GC=\{\gamma:\exists c\in C, c(\gamma)\neq 0\}$ is an open subgroupoid  of $G$ such that  $C=\spn N(B,C)= A_R(\GC;\Sigma|_\GC)$. So the map $C\mapsto \GC$ is the desired lattice isomorphism. 
\end{proof}
\section{Algebraic Cartan Pairs and Conditional Expectations}\label{CPCE}
Theorem~\ref{laticeIso} gives a correspondence between subalgebras $A_R(\go)\subseteq C\subseteq A_R(G)$ with  $A_R(\go)\subseteq C$ quasi-Cartan and wide open subgroupoids of $G$, but the question becomes: how can we tell if a given intermediate subalgebra has  $A_R(\go)\subseteq C$ quasi-Cartan? 
For algebraic Cartan pairs, one way to guarantee that the intermediate subalgebra also gives an ACP is if it is the image of a conditional expectation. 
We show this after the following lemma.

\begin{lem}\label{lem: nE(n) in commutant}
    Let $A$ be an $R$-algebra, and let $D \subseteq A$ be an abelian subalgebra with $I(D)$ being a set of local units for $A$.
    Suppose $C$ is an algebra with $D \subseteq C \subseteq A$ and assume that there is a conditional expectation $\exa\colon A \rightarrow C$.
    Then for $n \in N(A,D)$ we have that both $n\exa(n^\dag)$ and $\exa(n)n^\dag$ commute with every element of $D$.
\end{lem}

\begin{proof}
    Note that, since $I(D)$ are local units for $A$, $n^\dag n \in D$ for all $n \in N(A,D)$.
    Take $n \in N(A,D)$, and take any $d \in D$. Then
\begin{align*}
    d n\exa(n^\dag) &= d nn^\dag n \exa(n^{\dag}) \\
    &= nn^{\dag} d n \exa(n^{\dag}) \\
    &= n \exa( n^{\dag} d n n^{\dag}) \\
    &= n \exa( n^{\dag} n n^{\dag} d) \\
    &= n \exa(n^{\dag}) d.
\end{align*}
Hence $n \exa(n^{\dag})$ commutes with $D$.
A similar argument works for$n\exa(n^\dag)$.
\end{proof}

\begin{prop} \label{CEqC}Let $R$ be a commutative ring with identity and $(\Sigma,i, q)$ be a discrete twist over an effective ample Hausdorff groupoid $G$.  Suppose $C$ is an algebra such that $A_R(\go)\subseteq C\subseteq \Stein$.   If there exists a conditional expectation $\exa: \Stein\to C$, then $(C, A_R(\go))$ is an algebraic Cartan pair.
\end{prop}

\begin{proof}
The only part of the definition of an ACP that is not automatic is that the inclusion is regular, that is, that $\spn N(C,A_R(\go))=C$.  To show this, it suffices to show $\exa(N(\Stein, A_R(\go)))=N(C, A_R(\go))$.  Take $n \in N(\Stein,A_R(\go))$.
By Lemma~\ref{lem: nE(n) in commutant}, $n \exa(n^{\dag})$  and $\exa(n)n^\dag$ commute with $A_R(\go)$.
Since $G$ is effective the algebra $A_R(\go)$ is maximal abelian in $\Stein$. Hence $n\exa(n^{\dag})$ and $\exa(n)n^\dag$ are in $A_R(\go)$.
Thus, for all $n \in N(\Stein,A_R(\go))$
$$ \exa(n)\exa(n^{\dag}) = n\exa(n^{\dag}) = \exa(n)n^{\dag}\quad{\text{and}}\quad \exa(n^{\dag}) \exa(n)=\exa(n^{\dag}) n = n^{\dag}\exa(n).$$

Take $n \in N(\Stein,A_R(\go))$. Applying the above equalities to $n$ and $n^{\dag}$ we get
\begin{align*}
    \exa(n) &= \exa(nn^{\dag}n) = \exa(n) n^{\dag}n\\
    &= \exa(n)\exa(n^{\dag})n \\
    &= \exa(n) \exa(n^{\dag}) \exa(n),
\end{align*}
and similarly
$$ \exa(n^{\dag}) = \exa(n^{\dag}) \exa(n) \exa(n^{\dag}).$$
Now take any $d \in A_R(\go)$. 
We show that $\exa(n) d \exa(n^\dag) \in A_R(\go)$ for all $n \in N(\Stein,A_R(\go))$.
Take $n \in N(\Stein,A_R(\go))$
\begin{align*}
    \exa(n) d \exa(n^{\dag}) &= \exa(nd) \exa(n^{\dag}) \\
    &= \exa(ndn^{\dag}n) \exa(n^{\dag}) \\
    &= ndn^{\dag}\exa(n) \exa(n^{\dag}) \in A_R(\go).
\end{align*}
Thus $\exa(n)\in N(C, A_R(\go))$, with $\exa(n)^\dag = \exa(n^\dag)$.
Therefore, $\exa(N(\Stein, A_R(\go)))=N(C, A_R(\go))$ as desired.
\end{proof}

The following allows us to relate the existence of a conditional expectation $\exa: \Stein\to C$ to a property of the groupoid.  

\begin{prop}\label{CEC}
Let $R$ be a commutative ring with identity and let $(\Sigma,i, q)$ be a discrete twist over an ample Hausdorff groupoid $G$. Suppose $H\subseteq G$ is an open subgroupoid. Then there is a conditional expectation $\exa:A_R(G,\Sigma)\to A_R(H, \Sigma|_H)$ if and only if $H$ is closed.
\end{prop}

\begin{proof}
Recall that $A_R(H, \Sigma|_H)$ is a subalgebra of $A_R(G, \Sigma)$ by Lemma~\ref{lem: subgroupoid to subalgebra}.
 Assume there exists a conditional expectation $\exa:A_R(G,\Sigma)\to A_R(H, \Sigma|_H)$.  Pick a net $\{\gamma_i\}\subseteq H$ such that $\gamma_i\to \gamma\in G$.  We show that $\gamma\in H$.   Since $G$ is ample, there exists a compact open bisection $B \subseteq G$ such that $\gamma\in B$.   Without loss of generality, assume $\gamma_i \in B$ for all $i$.  Pick an $R$-linear homeomorphism $\psi_B: \lineP\inv(B)\to R\times B$ and consider $1_{\psi_B}$. 
For each $i$, there exists a compact open bisection $\gamma_i\in B_i\subseteq H$ and by intersecting with $B$, we can assume $B_i\subseteq B$ too.
For each $i$ we can assume $\psi_{B_i}=\psi_{B}|_{B_i}$ and so   $\chrp{r(B_i)} \chrp{B}= \chrp{B_i} \in A_R(H, \Sigma|_H)$.
Further $[\chrp{r(B_i)}h](\sigma) = h(\sigma)$ for all $\sigma \in B_i$ and $h \in \Stein$.
Hence
\[\exa(\chrp{B})(\gamma_i)=[\chrp{r(B_i)}\exa(\chrp{B})] (\gamma_i) = \exa(\chrp{r(B_i)}\chrp{B}) (\gamma_i) =  \chrp{B_i}(\gamma_i) = 1.\]
Since $\exa(\chrp{B})\in A(H,\Sigma|_H)$, by Proposition~\ref{prop regular} (\cite[Proposition~2.8]{A+8}), there exists a finite collection $\mathcal{F}$ of compact open bisections $\{U\}_{\mathcal{F}}$, $R$-linear homeomorphisms $\psi_U: R\times U\to \lineP\inv(U)$, and elements $\re_U\in R$ for each $U\in \cF$ with $\exa(\chrp{B})=\sum_{\mathcal{F}} \re_U \chrp{U}$.  So $\supp( \exa(\chrp{B}))\subseteq \bigcup_{\mathcal{F}} U\subseteq H$ is compact. Since $\{\gamma_n\}$ is a sequence in the compact set $\supp( \exa(\chrp{B}))\subseteq H$, it must have a convergent subnet $\gamma_{i_\ell}\to \eta\in H$.  We have that $G$ is Hausdorff and $\gamma_{i_\ell}\to \gamma$ by assumption and hence $\gamma=\eta\in H$ as well.  Thus $H$ is closed.

For the converse, suppose $H$ is closed.  Then   $\exa: A_R(G,\Sigma)\to A_R(H, \Sigma|_H)$ given by $\exa(f)=f|_{H}$ is the desired conditional expectation.
\end{proof}
\begin{corollary}\label{CE qC} Let $R$ be a commutative ring with identity and let $(\Sigma,i, q)$ be a discrete twist over an effective ample Hausdorff groupoid $G$. Then there exists an order preserving bijection between the algebras $C$ with $A_R(\go)\subseteq C\subseteq \Stein$, that have conditional expectations $\Stein\to C$ and wide clopen subgroupoids of $G$.
\end{corollary}

\begin{proof}
This is a consequence of Proposition~\ref{CEqC}, Theorem~\ref{laticeIso}, and Proposition~\ref{CEC}.
\end{proof}

The next example shows we can not extend Proposition~\ref{CEqC} (nor Corollary~\ref{CE qC}) to the quasi-Cartan case.  As both of these results are true for Cartan inclusions of $C^*$-algebras \cite[Theorem 3.5]{BEFPR},  these results illustrate both the differences between the theory of Cartan inclusions of $C^*$-algebras and the theory of quasi-Cartan inclusions of algebras and the distinction between Cartan inclusions and quasi-Cartan inclusions of algebras.

\begin{example}\label{ex: intermediate fail}
    Consider $A_\CC(\ZZ_3)$.  Then $\CC\subseteq A_\CC(\ZZ_3)$ is quasi-Cartan, but not Cartan.  Consider $C=\{\fe \delta_0+\fea(\delta_1+\delta_2):\fe,\fea\in \CC\}$. We show in Proposotion~\ref{bad apple} below that$\CC\subseteq C\subseteq A_\CC(\ZZ_3)$ is an intermediate subalgebra and the inclusion $\CC\subseteq C$ is not quasi-Cartan.  But there is a faithful conditional expectation $\exa: A_\CC(\ZZ_3)\to C$ given by 
    \[
    \fe_0\delta_0+\fe_1\delta_1+\fe_2\delta_2\mapsto \fe_0\delta_0+\frac{\fe_1+\fe_2}{2}\left(\delta_1+\delta_2\right).
    \]
    That $\exa$ is $\CC$-linear and $\exa|_C=\id$ are trivial.  That $\exa$ is $C$-linear and faithful are short computations.   
    Thus $C$ illustrates that the conclusions of Proposition~\ref{CEqC} can fail outside of the Cartan setting.
\end{example}
\section{Diagonal Pairs}\label{DP}
We showed in Section~\ref{SGC} that there is a lattice isomorphism between intermediate quasi-Cartan pairs and wide open subgroupoids.
In this section, we aim to answer the following question: 

\begin{ques}\label{ques main}
If $D\subseteq A$ is quasi-Cartan and $D\subseteq C\subseteq A$, must $D\subseteq C$ be quasi-Cartan?
\end{ques}

The next example shows that for this question to have a positive answer, our algebras must be over a field $K$, that is $R=K$ is a field.

\begin{example}\label{PQC implies field}   Let $D\subseteq A$ be a quasi-Cartan pair over a ring $R$.
Assume that $D \neq A$.
By \cite[Theorem~6.6]{A+8}, we can assume $A=\Stein$ and $D=A_R(\go)$.   Suppose that $R$ is not a field, that is, suppose there exists nonzero $\re\in R-R^\times$.   Let $C$ be the algebra generated by $D$ and $\{\re n: n\in N(A,D)\}$. 

Take any $\gamma \in G - \go$.
Then there exists $n \in N(A,D) - D$ such that $n(\gamma) = 1$.
Hence $\re n(\gamma) \neq 0$.
Therefore $\gamma \in G_C$, and hence $G_C = G$.
Further, it follows that $D \neq C$.
Note that $n \notin C$ for any $n\in N(A,D)-D$.
Hence, the inclusion $(C,D)$ is not a regular inclusion and so is not quasi-Cartan.

For a more explicit example, let $R=\ZZ$.  Since $\ZZ$ satisfies the local bisection hypothesis, $(A_\ZZ(\ZZ), A_\ZZ(\{0\}))$ is quasi-Cartan.  
Then let $A=A_\ZZ(\ZZ)\cong\ZZ[x,x\inv]$ and let $\re=2$.   Then $C$ as defined above is the algebra of Laurent polynomials where the coefficients of every nonzero power of $x$ are all even.
\end{example}

We now suppose $R=K$ is a field.  
Example~\ref{ex: intermediate fail} shows that the answer to Question~\ref{ques main}  is still no in general (see also Section~\ref{PqC}).  However, if we assume a stronger statement that $A$ is a $K$-algebra over a field $K$  and  $D\subseteq A$ is {\em diagonal},  we show in Theorem~\ref{thm: diagonal case}  that every intermediate subalgebra $D\subseteq C\subseteq A$, must have  $D\subseteq C$ diagonal.   
From Theorem~\ref{thm: twists aqp etc} (\cite[Proposition~7.1]{A+8}) we know that $D\subseteq A$ diagonal if and only if there exists a discrete twist $\Sigma\to G$ with  $G$ {\em principal}, such that $A\cong \Stein$ and $D\cong A_R(\go)$.
We will thus focus on twisted Steinberg algebras over principal groupoids in the remainder of this section.  To prove our main theorem of this section we need a lemma.
\begin{lem} \label{normalizers in C} Let $K$ be a field and $\Sigma\to G$ be a discrete twist over a principal ample Hausdorff groupoid $G$.   Suppose that $c\in \Steink$ and $\gamma\in G$ such that $c(\gamma)\neq 0$.   
Then there exists a normalizer $n$ in the $A_K(\go)$-bimodule generated by $c$ (denoted $\bi(c)$) such that $c|_{S_n}=\re n$ for some $\re\in K-\{0\}$.  Moreover,
\[c\in \spn_K N(\Steink, A_K(\go))\cap \bi(c).\]
\end{lem}

\begin{proof}
By Proposition~\ref{prop regular}(\cite[Corollary~2.10]{A+8}) there exists disjoint compact open bisections $B_i$, $K$-linear homeomorphisms $\psi_{B_i}: K\times B_i\to \lineP\inv(B_i)$ and $\fe_i\in K$ such that  $ c=\sum \fe_i \chrp{B_i}$.  Thus, since $c(\gamma)\neq 0$, we may assume there exists a unique $i_0$ such that $\gamma\in B_{i_0}$.  Let $\{V_k\}$ and $\{W_k\}$ be compact open neighborhood bases of $r(\gamma)$ and $s(\gamma)$ respectively.  We claim that there exists $N_i$ such that for $k>N_i$, $V_kB_i W_k=\emptyset$ for all $i\neq i_0$.  Consider $B_i s(\gamma)$. There are two cases.

Case 1: $B_i s(\gamma)=\emptyset$. Then there exists $N_i$ such that $B_i W_k=\emptyset$ for all $k>N_i$ as $B_i$ is compact.  

Case 2: $B_is(\gamma)$  is nonempty.  Then  there exists $\eta_i$ such that $\{\eta_i\}=B_i s(\gamma)$ as $B_i$ is a bisection.  But $G$ is principal so $r(\eta_i)\neq r(\gamma)$, and hence there exist disjoint open sets $R_i, S_i$ such that $r(\gamma)\in R_i$ and $r(\eta_i)\in S_i$.  Each of $\{V_k\}$ and $\{W_k\}$ are neighborhood bases, so there exists $N_i$ such that for $k>N_i$, $V_k\subseteq R_i$ and $W_k\subseteq s(S_iB_i)$.  Thus $V_k B_i W_k=\emptyset$ for $k>N_i$. 

Take $N=\max\{N_i\}$, pick $k>N$ and take $V=V_k$ and $W=W_k$. Then  
\[1_V c 1_W=\sum \fe_i 1_V \chrp{B_i} 1_W=\fe_{i_0}  1_V \chrp{B_{i_0}} 1_W.\]
Now take $n:=1_V \chrp{B_{i_0}} 1_W$  and $\fe=\fe_{i_0}$. Then $n$ a normalizer and $n=\fe_{i_0}\inv 1_V c1_W$ is in the $A_K(\go)$-bimodule generated by $c$.  We have  $\gamma \in VB_{i_0}W$  and $c|_{S_n}=\fe n$ by construction.

To show the last statement, for each $\gamma\in S_c$ find $V_\gamma, W_\gamma$  so that $1_{V_{\gamma_k}}c1_{W_{\gamma_k}}$ is a scalar multiple of a normalizer in $\bi(c)$ as in the previous paragraph.  The support of $c$ is compact  and $\{V_\gamma B_i W_\gamma\}_{\gamma, i}$  forms an (compact) open cover of this support.  So there is a finite subcover $\{V_{\gamma_k} B_{i_k} W_{\gamma_k}\}$.    We now have 
\[
c=\sum_k 1_{V_{\gamma_k}}c1_{W_{\gamma_k}} \in  \spn_K N(\Steink, A_K(\go))\cap \bi(c)).\qedhere
\]
\end{proof}

We now present the main theorem of this section.
\begin{thm}\label{thm: diagonal case} Let $K$ be a field  and let $\Sigma\stackrel{q}{\to} G$ be a discrete twist over a principal ample Hausdorff groupoid $G$ such that $(A_K(G; \Sigma),A_K(\go))$ is a diagonal inclusion.   Suppose $C$ is an intermediate subalgebra, $A_K(\go)\subseteq C\subseteq \Steink$.  Then $A_K(\go)$ is diagonal in $C$. \end{thm}

\begin{proof}
We show that $(C,A_K(\go))$ is a quasi-Cartan inclusion. By Lemma~\ref{groupoid inclusions}, it suffices to show $C=\NC$. That $\NC\subseteq C$ is immediate. The reverse inclusion follows from Lemma~\ref{normalizers in C}: for  $c\in C$,  $\bi(c)\subseteq C$ and so  $N(\Steink, A_K(\go))\cap \bi(c)\subseteq N(C,A_K(\go))$, and thus $c\in \spn N(C,A_K(\go))=\NC$.

 By Theorem~\ref{laticeIso}  we have $A_K(\Sigma|_\GC;\GC)=C$ where $\GC=\{\gamma\in G: \exists c\in C, c(\gamma)\neq 0\}$ and $\Sigma|_\GC=q\inv(\GC)$.    Since $G$ is principal, $G_C$ is principal and we get $A_K(\go)\subseteq C=A_K( \Sigma|_\GC;\GC) $ is diagonal from Corollary~\ref{passes to open}.   
\end{proof}

\begin{rmk} Theorem~\ref{thm: diagonal case} will be subsumed by Theorem~\ref{i2i cartan} below, where we extend the groupoids to which the theorem applies. We include Theorem~\ref{thm: diagonal case} here for several reasons. \begin{itemize}
    \item Theorem~\ref{thm: diagonal case} is the direct algebraic analogue of the existing theorem in the $C^*$-algebraic literature \cite[Theorem~4.5]{BEFPR}.
    \item Theorem~\ref{thm: diagonal case} serves as motivation for our results in Section~\ref{PqC} below.
    \item Theorem~\ref{thm: diagonal case} relates to our results on $A_R(\go)$-bimodules in Section~\ref{bimodules} where we show in Corollary~\ref{cor: biiff} $A_R(\go)$-bimodules correspond to open sets of the groupoid if and only if $G$ is principal.  Indeed Lemma~\ref{normalizers in C} is used for both results.
\end{itemize}
\end{rmk}
\section{Purely quasi-Cartan}\label{PqC}
\begin{defn} We say an inclusion $D\subseteq A$ is {\em  purely quasi-Cartan} if for every intermediate subalgebra $D\subseteq C\subseteq A$ we have $D\subseteq C$ is quasi-Cartan.  
\end{defn}

Theorem~\ref{thm: diagonal case} shows that $A_K(\go)\subseteq \Steink $  is purely quasi-Cartan if $G$ is principal (or equivalently diagonal inclusions are always purely quasi-Cartan).
In this section we explore when quasi-Cartan inclusions are purely quasi-Cartan.  It turns out that there are many examples of purely quasi-Cartan inclusions, but that this class is still quite restrictive.   The next proposition shows that $(A_R(G; \Sigma),A_R(\go))$ can be purely quasi-Cartan only if any pair of distinct  elements in the unit space have at most one groupoid element connecting them. 

\begin{prop}\label{prop: two arrows} Suppose that $R$ is a commutative unital ring and suppose $\Sigma\to G$ is a discrete twist over an ample Hausdorff groupoid $G$ such that there exists $u,v\in \go$ with $u\neq v$ and  $|uGv|\geq 2$.  Then there exists an intermediate subalgebra $A_R(\go)\subseteq C\subseteq \Stein$ such that $(C,A_R(\go))$ is not quasi-Cartan. 
\end{prop}

\begin{proof}
Let $\gamma_1,\gamma_2\in uGv$ with $\gamma_1\neq \gamma_2$ and pick $B_1, B_2$ disjoint open bisections containing $\gamma_1, \gamma_2$ respectively and fix  $U_u, U_v$ disjoint open sets containing $u,v$ respectively.  Take $\Ooo_i=U_uB_i U_v$.   Pick $R$-linear homeomorphisms $\psi_{\Ooo_i}: R\times \Ooo_i\to \lineP\inv(\Ooo_i)$. Define
\[
f=\chrp{\Ooo_1}+\chrp{\Ooo_2}.  
\]
Using that $r(\Ooo_i)\cap s(\Ooo_i)=\emptyset$ we have  $f^2\equiv 0$.   Define $C$ to be the algebra generated by $A_R(\go)$ and $f$.  Then 
\[
C=\spn\{g_1 fg_2: g_i\in A_R(\go)\}.
\]
We claim $A_R(\go)\subseteq C$ is not quasi-Cartan.  To see this we show that for all $h\in C$, $h(\gamma_1)=h(\gamma_2)$ so that $C$ is not given by $A(H)$ for any wide open subgroupoid $H$ of $G$.  This is a quick computation as 
\[
(g_1fg_2)(\gamma_1)=g_1(u)f(\gamma_1)g_2(v)=g_1(u)f(\gamma_2)g_2(v)=(g_1fg_2)(\gamma_2).
\]  The result now follows by applying Theorem~\ref{laticeIso}.  
\end{proof}

\begin{rmk}
If $G$ is not principal, then there exists $v$ such that $|vGv|\geq 2$.  If $G$ is also strongly effective there must exists $u$ such that $uGv\neq \emptyset$, otherwise $v$ would be a closed invariant set consisting of isotropy.  Pick $\gamma\in uGv$ and $\eta_1\neq \eta_2\in vGv$ then $\gamma\eta_1\neq\gamma\eta_2\in uGv$ so that $|uGv|\geq 2$.  Proposition~\ref{prop: two arrows} shows that non-principal strongly effective groupoids do not have a purely quasi-Cartan Steinberg algebras. \end{rmk}

Next we show groupoids with large isolated isotropy groups do not yield purely quasi-Cartan inclusions either.  
 We obtain this as a corollary to the following two results. Before stating the next lemma observe that if $v\in \go$ with $vG=Gv=vGv$, then $\{v\}$ is an invariant set and so Lemma~\ref{lem:exact} says  that restriction of functions gives a homomorphism $p:\Stein\to A_R(vGv;v\Sigma v)$ of $R$-algebras. 

\begin{lem}\label{Bad apple} Let $R$ be a commutative unital ring and  $(\Sigma, i, q)$ be a discrete twist over an ample Hausdorff groupoid $G$  and let $\Gamma$ be a subgroup of $vGv$ for some $v\in \go$ with $vG=Gv=vGv$.  Let $p:\Stein\to A_R(vGv;v\Sigma v)$ be the homomorphism of $R$-algebras obtained by restriction of functions. Suppose that there exists a subalgebra $A\subsetneq A_R(\Gamma; \Sigma|_\Gamma)$ such that 
\begin{enumerate}
\item $A$ contains the identity of $A_R(\Gamma; \Sigma|_\Gamma)$,
\item for every $s\in \Gamma$ there exists  $b\in A$ such that $b(s)\neq 0$,
\item there exists $s_0\in \Gamma$ such that  $\delta_{s_0}\notin A$.
\end{enumerate}
If we take $C:=p\inv(A)$, then $A_R(\go)\subseteq C\subsetneq A_R(\GC;\Sigma|_\GC)=p\inv(A_R(\Gamma; \Sigma|_\Gamma))$: in particular $A_R(\go) \subseteq C$ is not quasi-Cartan.
\end{lem}

\begin{proof} First note that since $p:\Stein\to A_R(vGv;v\Sigma v)$ is a homomorphism, $C$ is a subalgebra of $\Stein$. Let $\Omega=\go-\{v\}$, then $\Omega$ is invariant.  Consider \[G\Omega=s\inv(\go-\{v\})=r\inv (\go-\{v\}).\] We start by showing $\GC=G\Omega\cup \Gamma$.   First, suppose that $\gamma\in G\Omega$. Pick a compact open bisection $B$ containing $\gamma$ such that $v\notin r(B)\cup s(B)$.   Then for any $R$-linear homeomorphism $\psi_B: R\times B\to \lineP\inv(B) $ we have $p(\chrp{B})=0\in A$ so that $\chrp{B}\in C$ and  since $\chrp{B}(\gamma)\neq 0$ we have $\gamma\in \GC$.    

Now suppose instead that $\gamma\in vGv$.  If $\gamma\in \Gamma$, then there exists $b\in A$ such that $b(\gamma)\neq 0$ and we can write $b=\sum_{\eta\in F\subseteq \Gamma} \re_\eta \delta_\eta$ for some finite subset $F \subseteq \Gamma$ with $\gamma\in F$.  Choose disjoint compact open bisections $B_\eta$ in $G$ around each $\eta\in F$.    Then there exists $R$-linear homeomorphisms $\psi_{B_\eta}: R\times B_\eta\to \lineP\inv(B_\eta)$ and $\re_\eta\in R$ such that 
\[\tilde{b}=\sum_{\eta\in F\subseteq \Gamma} \re_\eta \chrp{B_\eta}\in q\inv (b)\]
and we have $\tilde{b}(\gamma)=\re_\gamma \chrp{B_\gamma}\neq 0$.
Hence $\gamma \in \GC$

Finally, we show that if $\gamma\in vGv-\Gamma$, then $\gamma\notin \GC$.   But if $a\in \Stein$ with $a(\gamma)\neq 0$, then $p(a)(\gamma)\neq 0$ so $p(a)\notin A_R(\Gamma; \Sigma|_\Gamma)$ and so $a\notin C$ and hence $\gamma\notin \GC$.      This proves that $\GC=G\Omega\cup \Gamma$ and $A_R(\GC;\Sigma|_{\GC})=p\inv(A_R(\Gamma))$. 

 To see that the containment  $C\subseteq A_R(\GC;\Sigma|_{\GC})$ is proper, pick  $B \subseteq G$  a compact open bisection such that $s_0\in B$ and  an $R$-linear homeomorphism $\psi_B: R\times B\to \lineP\inv(B) $.  Then $p(\chrp{B})=\delta_{s_0}\notin A$, so that $\chrp{B}\notin C$.  Since $\chrp{B}\in A_R(\GC;\Sigma|_{\GC})$ we must have $C\subsetneq A_R(\GC;\Sigma|_{\GC})$.  
\end{proof}

\begin{prop}\label{bad apple} Let $R$ be a commutative ring with identity such that $1\neq -1$. Let $(\Sigma, i, q)$ be a discrete twist over an ample Hausdorff groupoid $G$ and let $v\in \go$ be such that $vG=Gv=vGv$ with $3\leq n=|vGv|$ and $v\Sigma v\to vGv$ is trivial.  Then there exists a subalgebra $A$ of $A_R(vGv)$ such that 
\begin{enumerate}
\item $A$ contains the identity of $A_R(vGv)$,
\item for every $s\in vGv$ there exists  $b\in A$ such that $b(s)\neq 0$,
\item and there exists $s_0\in vGv$ such that  $\delta_{s_0}\notin A$.
\end{enumerate}
\end{prop}

\begin{proof} 
Let $v\in \go$ as in the statement and pick $\gamma\in vGv-\{v\}$.  Take $\Gamma$ to be the cyclic group generated by $\gamma$.   There are two cases.

Case 1: $|\Gamma|=\infty$.  Then $A_R(\Gamma)\cong A_R(\ZZ)\cong R[x]$.  Let $\tilde A$ be the subalgebra of $R[x]$ generated by all $f\in R[x]$ such that $f( -1)=f(1)$.  Note that if $i$ is even we have $1^i=(-1)^i$ and if $i$ is odd then $1^i-1^{-i}=0=(-1)^i-(-1)^{-i}$.  Pulling $\tilde{A}$ back to a subalgebra $A$ of $A_R(\ZZ)$ we get that $\delta_i\in A$ if $i$ is even, $\delta_i-\delta_{-i} \in A$ if $i$ is odd, and $\delta_i\notin A$ if $i$ is odd.

Case 2: $|\Gamma|=n< \infty$.  Then $\Gamma\cong \mathbb{Z}_n$.  Let $A=\left\{\re_v\delta_v+\re\left(\sum_{i\neq v} \delta_i\right)\right\}$.  Then we claim that $A$ is the desired subalgebra.   First, we show that $A$ is a subalgebra
\begin{itemize}
\item $A$ is closed under sums $\re_v\delta_v+\re\left(\sum_{i\neq v} \delta_i\right)+\re'_v\delta_v+\re'\left(\sum_{i\neq v} \delta_i\right)=(\re_v+\re'_ v)\delta_v +(\re+\re')\left(\sum_{i\neq v} \delta_i\right)$
\item $A$ is closed under scalar products: $\tilde{\re}(\re_v\delta_v+\re\left(\sum_{i\neq v} \delta_i\right))=\tilde{\re}\re_v\delta_v+\tilde{\re}\re\left(\sum_{i\neq v} \delta_i\right)$
\item $A$ is closed under products:
\begin{align*}
\bigg( \re_v\delta_v&\left.+\re\left(\sum_{i\neq v} \delta_i\right)  \right) \left(\re'_v\delta_v+\re'\left(\sum_{j\neq v} \delta_j\right)\right)\\
&=\re_v\re'_v\delta_v+\re_v\re'\left(\sum_{i\neq v} \delta_i\right)+\re_v'\re\left(\sum_{i\neq v} \delta_i\right)+\re\re'\left(\sum_{i\neq v} \delta_i\right)\left(\sum_{j\neq v} \delta_j\right)\\
&=\re_v\re'_v\delta_v+(\re_v\re'+\re_v'\re)\left(\sum_{i\neq v} \delta_i\right)+\re\re'\left(\sum_{i\neq v, j\neq v} \delta_{ij}\right)\\
&=\re_v\re'_v\delta_v+(\re_v\re'+\re_v'\re)\left(\sum_{i\neq v} \delta_i\right)+\re\re'\left(\sum_{k}\sum_{j\neq v} \delta_{kj\inv j}\right)\\
&=\re_v\re'_v\delta_v+(\re_v\re'+\re_v'\re)\left(\sum_{i\neq v} \delta_i\right)+\re\re'\left(\sum_{k}(n-1)\delta_{k}\right)\\
&=(\re_v\re'_v+\re\re'(n-1))\delta_v +(\re_v\re'+\re_v'\re+\re\re'(n-1))\left(\sum_{k\neq v}\delta_{k}\right).
\end{align*}
\end{itemize}
Thus we get $A$ is an algebra.  Taking $\re_v=1$ and $\re=0$ we see that $A$ contains the identity of $\Gamma$.   Letting $\re=\re_v=1$ we get $b=\sum_{k} \delta_k\in A$ and $b(k)=1$ for all $k\in \Gamma$.     Lastly pick $s_0\neq v$.  We claim $\delta_{s_0}\notin A$. Pick $s_1\in vGv-\{s_0,v\}$ such an $s_1$ exists since $|vGv|\geq 3$.   Then $\delta_{s_0}(s_0)=1\neq 0=\delta_{s_0}(s_1)$ but by definition every element of $A$ has the same value on all elements of $vGv-\{v\}$.   Thus $\delta_{s_0}\notin A$ as desired.
\end{proof}

\begin{corollary}\label{iso small}Let $R$ be a commutative ring with identity such that $1\neq -1$. Suppose $G$ is an ample Hausdorff groupoid and $\Sigma=R^\times \times G\to G$  is a trivial twist such that  for some $v\in \unit{G}$ we have $vG=Gv=vGv$ with $3\leq n=|vGv|$. Then there exists a subalgebra $C$ with $A_R(\go)\subseteq C\subseteq A_R(G)$ and $A_R(\go)\subseteq C$ not quasi-Cartan.
\end{corollary}

\begin{proof}Apply Lemma~\ref{Bad apple} to the algebra obtained in Proposition~\ref{bad apple}.
\end{proof}

\begin{example}
Let $R$ be a commutative ring  with identity  such that  $1\neq -1$.  Consider the graph $E$, drawn below. 
\begin{figure}  [h!]
\xymatrix{
\bullet_v \ar@(ul,ur)^{e_1}\ar@(ul,dl)_{e_2}\ar@(dr,dl)^{e_3}\ar[r]_f & \bullet_w\ar@(ur,dr)^g
}
\end{figure}
Let $A_R(G_E)$ be the Steinberg algebra associated to $G_E$ over $R$. 
Since every loop has an entrance, $E$ satisfies condition $L$, so the associated groupoid $G_E$ is topologically principal.   Therefore, by Theorem~\ref{thm: twists aqp etc} (\cite[Proposition 7.1]{A+8}, $A_R(G_E^{(0)})$ is Cartan in $A_R(G_E)$.

Now if we take $x=gggggg...$, then $xG=Gx=xGx$ and $|xGx|=\infty\geq 3$ and so Corollary~\ref{iso small} shows that $A_R({G_E}^{(0)})\subseteq A_R(G_E)$ is not purely quasi-Cartan.
\end{example}

Combining Corollary~\ref{iso small} and Proposition~\ref{prop: two arrows}, we see that for $A_R(\go)\subseteq \Stein$ to be purely quasi-Cartan the following must be satisfied:  

\begin{equation}\label{I2I} \text{If $u,v \in \go$ such that $|vGw|>1$, then $v=w$, $vG=Gv=vGv$ and $|vGv|=2$. }\end{equation}  
\begin{defn}\label{def:I2I}We call a Hausdorff \'etale groupoid satisfying \eqref{I2I}  an  \emph{$\i2i$-groupoid} (Isolated $\ZZ_2$-Isotropy). 
\end{defn}

\begin{example} \label{Switch} For an explicit example of a non-principal $\i2i$-groupoid, let $\ZZ/2\ZZ$ act on $[-1,1]$ by switching signs.    Let $Y$ consist of the cantor set $\mathcal{C}$ union $\mathcal{C}-1$ and restrict the $\ZZ/2\ZZ$ action to $Y$.  
Let $G$ be the associated groupoid, then $G$ is a non-principal ample $\i2i$-groupoid.
\end{example}

We will show that essentially all purely quasi-Cartan inclusions come from $\i2i$-groupoids. 

\begin{lem}\label{i2i normalizers} Let $K$ be a field and $\Sigma\to G$ be a discrete twist over an ample $\i2i$-groupoid $G$.  Suppose  $c\in\Steink$ and $\gamma\in G$ such that $c(\gamma)\neq 0$. Then there exists a normalizer $n$ in the subalgebra $\su(c)$  generated by $A_K(\go)$ and $c$, such that $c|_{S_n}=\fe n$ for some $\fe\in K-\{0\}$. Moreover $c\in \spn_K N(\Steink, A_K(\go))\cap \su(c)$. \end{lem}

\begin{proof} Since $c\in \Steink$ we can assume $c=\sum_{B\in F} \fe_B \chrp{B}$ where $B$ are all disjoint, $F$ is finite,  $\psi_B: K\times B\to \lineP\inv(B)$ are $K$-linear homeomorphisms, $\fe_B\in K$.  Since $\su(c)$ contains $A_K(\go)$, by subtracting elements of $A_K(\go)$ we can assume $\supp c\subseteq G-\go$.   For $\gamma\in \supp c$, we show there exists $V,W\subseteq \go$, $B_0\in F$ such that  $\gamma\in VB_0W$ and $1_V c1_W=\fe_{B_0} \chr{\psi_{B_0}|_{VB_0W}}$.  Consider $\gamma\in \supp c$.  There are two cases.

Case 1:   $r(\gamma)\neq s(\gamma)$.   There are at most finitely many $\gamma_i$ such that $\gamma_i\in \supp(c)$ and $s(\gamma_i)=s(\gamma)$.  Define $B_i$ to be the element of $F$ such that $\gamma_i\in B_i$ and $B_0$ to be the element of $F$ such that $\gamma\in B_0$.  By assumption $r(\gamma_i)\neq r(\gamma_j)$ for $i\neq j$ and $r(\gamma_i)\neq s(\gamma_i)$ for all $i$.  Thus we can pick mutually disjoint compact open neighborhoods $U_i$ of $r(\gamma_i)$ and $W_0$ of $s(\gamma)$ such that $U_i\cap W_0=\emptyset$.   Now let $W_i=s(U_i B_i)$ for $i\neq 0$.  Take $W=\cap_{i} W_ i$ and $V=U_0$.  Then $1_V c1_W=1_V (\fe_{B_0} \chrp{B_0})1_W=\fe_{B_0} \chr{\psi|_{B_0}|_{VB_0W}}$.  

Case 2: $r(\gamma)=s(\gamma)$.  Then by assumption $r(\gamma)G=Gr(\gamma)=r(\gamma)Gr(\gamma)=\{r(\gamma),\gamma\}$.   Thus there exists a unique $B_0\in F$ such that $r(\gamma)\in r(B_0)$.  Notice $B_0$ is also the only element of $F$ so that $s(\gamma)\in s(B_0)$.   Now since $F$ is finite and each $B$ is compact open and that $\go$ is Hausdorff we have that $\Ooo=\bigcup_{B\in F-\{B_0\}} (r(B)\cup s(B))$ is closed compact and open.  Then $V=r(B_0)-\Ooo$ and $W=s(B_0)-\Ooo$ are compact open.  Now $1_{V} f1_{W}=1_{V} (\fe_{B_0} \chrp{B_0})1_{W}=\fe_{B_0} \chr{\psi_{B_0}|_{VB_0W}} $ and $f(\gamma)=\fe_{B_0} \chr{\psi_{B_0}|_{VB_0W}}(\gamma).$

In either case we have $\gamma\in VB_0W$ and $1_V c1_W=\fe_{B_0} \chr{\psi_{B_0}|_{VB_0W}}$.  Since $1_V, 1_W\in A_K(\go)\subseteq \su(c)$ we have $n:=\chr{\psi|_{B_0}|_{VB_0W}}=\fe_{B_0}\inv 1_Vc1_W$  is a normalizer in $\su(c).$  By construction $c|_{S_n}=\fe_{B_0}n$ giving the result.
\end{proof}

\begin{thm}\label{i2i cartan}Let $K$ be a field and  let $(\Sigma, i, q)$ be a discrete twist over an ample $\i2i$-groupoid $G$. Then for every subalgebra $A_R(\go)\subseteq C\subseteq \Steink$ we have that $A_R(\go)\subseteq C$ is Cartan and so $A_R(\go)\subseteq \Steink$ is purely quasi-Cartan.
\end{thm}

\begin{proof}Follow the proof of Theorem~\ref{thm: diagonal case}  replacing the use of Lemma~\ref{normalizers in C} with Lemma~\ref{i2i normalizers}.
\end{proof}

\begin{corollary}\label{cor i2i trivial}
Suppose $R$ is a commutative ring with identity such that $1\neq -1$, $G$ a Hausdorff ample groupoid and $R^\times\times G\to G$ is a trivial twist.   Then $A_R(\go)\subseteq A_R(G)$ is purely quasi-Cartan if and only if $G$ is an $\i2i$-groupoid and $K$ is a field.
\end{corollary}

\begin{proof}
The ``if'' direction is Theorem~\ref{i2i cartan}.  For the ``only if'' direction combine Corollary~\ref{iso small}, Proposition~\ref{prop: two arrows}, and Example~\ref{PQC implies field}.
\end{proof}

\section{Bimodules}\label{bimodules}

In this section we look more generally to characterize $A_R(\go)$-bimodules in $\Stein$.   For an open subset $U \subseteq G$, we denote
\[A_R(G;\Sigma; U)=\spn\{\chrp{B}:B\subseteq U\}.\]
We can identify those $A_R(\go)$-bimodules of the form $A_R(\Sigma;G; U)$.

\begin{prop}\label{prop: left modules}
Let $R$ be a commutative ring with identity and let $(\Sigma,i,q)$ be a discrete twist over an ample Hausdorff groupoid $G$.  Suppose for  $D=A_R(G^{(0)})$ and $A=\Stein$ that $(A,D)$ is an AQP, in particular  $D$ is  $\WT$ with respect to $R$.
Let $M$ be a left $A(G^{(0)})$-module.
Then $M = A_R(G;\Sigma;U)$ (and hence  is a bimodule) for an open subset $U$ if and only if $M$ is the span of the normalizers in $N(\Stein, A_R(G^{0})) \cap M$.
\end{prop}

\begin{proof}
If $M = A_R(G; \Sigma;U)$ for some open set $U \subseteq \go$, that $M$ is a left $A(G^{(0)})$-module is straightforward. Conversely, suppose that $M$ is the span of the normalizers in $N(\Stein, A_R(G^{0})) \cap M$ and let
$$ U = \bigcup_{n \in N(\Stein,A(G^{0})) \cap M}S_n $$
where $S_n\subseteq G$ is the support of $n$. Then $M \subseteq A_R(G; \Sigma;U)$.
Moreover, by assumption and Proposition~8.6~\cite{A+8} we have, 
\[ M = \spn\{\chrp{B} \colon \chrp{B} \in M\}. \]

Take any compact open bisection $B \subseteq U$.
As $B \subseteq U$ and $B$ is compact, there are compact open bisections $C_1,\ldots,C_k$  and $R$-linear homeomorphisms $\psi_{C_i}: R\times C_i\to \lineP\inv (C_i)$ such that $B \subseteq \bigcup_{i=1}^k C_i$ and $\chrp{C_i} \in M$ for $1\leq i \leq k.$
Note that, since $M$ is a left $A(G^{(0)})$-module, we have that if $\chrp{Y} \in  M$ for an compact open bisection $Y \subseteq G$, and $X\subseteq Y$ is an compact open bisection, then $1_{{\psi_Y}|_{X} }\in M$.
Indeed, $1_{\psi_Y|_{X} }= 1_{r(X)}\chrp{Y}\in M$.

Hence $ 1_{\psi_{C_i}|_{B}}\in M$ for $1\leq i \leq k$.  Now for each $\gamma\in B$ there exists $C_{i(\gamma)}$ such that $\gamma\in C_{i(\gamma)}\cap B$ and by Lemma~\ref{lem: characteristic} there exists $\Ooo_\gamma\subseteq C_i\cap B$ and $\re_\gamma$ such that $\re_\gamma \chrp{C_i}|_{\Ooo_\gamma}=\chrp{B}|_{\Ooo_\gamma}$.  From above we get $\re_\gamma \chrp{C_{i(\gamma)}}|_{\Ooo_\gamma}\in M$.  Since $B$ is compact we get a finite collection of bisections $\{\Ooo_{\gamma_j}\}$ and elements of $R$, $\{\re_{\gamma_j}\}$ so that 
\[ 1_B = \sum \re_{\gamma_j}\chrp{C_{i(\gamma_j)}}|_{\Ooo_{\gamma_j}} \in M. \]
Thus $A_R(G; \Sigma;U) \subseteq M$, completing the proof.
\end{proof}

\begin{thm}[Spectral Theorem for Bimodules]\label{STbi}
Suppose $K$ is a field and $\Sigma\to G$ is a discrete twist over a principal ample Hausdorff groupoid $G$. Then the map $U \mapsto A_K(G;\Sigma; U)$ gives a one-to-one correspondence from open subsets of $G$ and $A_K(G^{(0)})$-bimodules in $\Steink$.
\end{thm}

\begin{proof}
Injectivity is straightforward. For surjectivity, let $M$ be an $A_K(\go)$-bimodule in $\Steink$.
Take $m \in M$, then $m\in \spn N(\Steink, A_K(\go))\cap \bi(m)$ by Lemma~\ref{normalizers in C}.
Hence $M$ is a span of normalizers in $N(\Steink, A_K(\go))$.
It follows from Proposition~\ref{prop: left modules} that $M = A_K(G;\Sigma;U)$ for an open subset $U \subseteq G$.
\end{proof}

Note that Theorem~\ref{thm: diagonal case} is a special case of Theorem~\ref{STbi}.
In Theorem~\ref{i2i cartan} we showed that we could extend Theorem~\ref{thm: diagonal case} to $\i2i$-groupoids.  We show below that the spectral synthesis given by Theorem~~\ref{STbi} works only in the diagonal case.

\begin{prop}[Failure of Spectral Synthesis]\label{prop: synth fail}
   Let $K$ be a field.  Suppose $G$ is an ample Hausdorff groupoid that  not principal 
    and $(\Sigma, i, q)$ is a discrete twist over $G$.
    There is a $A_K(\go)$-bimodule $M \subseteq \Steink$ such that $M$ is not equal to $A_K(G;\Sigma;U)$ for any open $U \subseteq G$.
\end{prop}

\begin{proof}
Chose $\gamma \in G$ such that $\gamma \notin \go$ and $r(\gamma) = s(\gamma)$.
Such a $\gamma$ exists, by the assumption that $G$ is not principal.
Choose disjoint compact open bisections $B_0$ and $B_1$ such that $r(\gamma) \in B_0$ and $\gamma \in B_1$ and  $K$-linear homeomorphisms $\psi_{B_i}: K\times B_i\to \lineP\inv(B_i)$.
Let $c = \chrp{B_0} + \chrp{B_1}$.
Let $\bi(c)$ be the $A_K(\go)$-bimodule generated by $c$.
Then
$$ \bi(c) = \spn\{1_V c 1_W \colon V,W \subseteq \go \text{ compact open}\}.$$
Let $V$ and $W$ be compact open subsets of $\go$ and let $a = 1_V c 1_W$ and let $S_a \subseteq G$ be the support of $a$.
Then $r(\gamma) \in S_a$ if and only if $\gamma \in S_a$ if and only if $r(\gamma) \in V$ and $r(\gamma) \in W$.
Thus, for any $a \in \bi(c)$, $a(r(\gamma)) = a(\gamma)$.
Hence, there is no open $W \subseteq G$ such that $\bi(c) = A_K(G;\Sigma;U).$
\end{proof}


\begin{corollary}\label{cor: biiff}
    Let $K$ be a field.  Suppose $\Sigma\to G$ is a discrete twist over an ample Hausdorff groupoid $G$. Then the map $U \mapsto A_K(G;\Sigma; U)$ gives a one-to-one correspondence from open subsets of $G$ and $A_K(G^{(0)})$-bimodules in $\Steink$ if and only if $G$ is principal.
\end{corollary}

\begin{proof}One direction is Theorem~\ref{STbi}, the other direction is Proposition~\ref{prop: synth fail}.
\end{proof}


\section{Intermediate subalgebras in $C^*$-Cartan inclusions}\label{cstar}

Quasi-Cartan inclusions were developed as algebraic generalizations of the notion of Cartan inclusions of $C^*$-algebras (which in turn was inspired by the von Neumann algebra setting).  
In \cite{B+6}, the authors show that in a sense, Cartan and quasi-Cartan are equivalent in the $C^*$-setting.   In this section, we assume familiarity with  basic $C^*$-algebraic theory as outlined, for example, in \cite{Murphy}.  

\begin{defn}[{\cite[Definition~5.1]{Ren08}}]\label{def: C*-Cartan} Suppose $A$ is a $C^*$-algebra and $D$ is a subalgebra of $A$,  we say $D$ is {\em Cartan} in $A$, or the pair $(A,D)$ is Cartan, if 
\begin{enumerate}
\item $D$ contains an approximate unit for $A$;
\item $A=\overline{\spn} N(A,D)$ where $N(A,D)=\{n\in A: nDn^*\cup n^*Dn\subseteq D\}$;
\item there exists a faithful conditional expectation $\ex: A\to D$;
\item $D$ is a maximal abelian subalgebra of $A$.
\end{enumerate}
\end{defn}

Note the definition of normalizers in the $C^*$-setting is different than the algebraic definition: $n^*n$ will be in $D$ but is unlikely to be an idempotent.  In particular, $n^{\dagger}$ and $n^*$ are different in general.
In the analogous definition for von Neumann algebras, one does not consider all normalizers, but only the normalizers which are partial isometries, see e.g. \cite[Definition~3.1 and Remark~3.2]{DFP}. 
When $n$ is a normalizer and a partial isometry, then $n^*n$ is a projection and $n^* = n^\dag$.

\begin{defn} We say an inclusion $(A,D)$ is {\em purely Cartan} if for every intermediate subalgebra $D\subseteq C\subseteq A$ we have $D\subseteq C$ is Cartan in $C$.\end{defn}

Using this definition,  \cite[Theorem~4.5]{BEFPR} states that if $(A,D)$ is diagonal and $A$ is nuclear then $(A,D)$ is purely Cartan.  Inspired by our progress on purely quasi-Cartan inclusions in the algebraic setting, we revisited \cite{BEFPR} to see what more can be said about purely Cartan inclusions in the $C^*$-setting.  We first outline the construction of Cartan pairs from twists.

A twist is a central extension
\[
\TT\times \go \stackrel{i}{\to} \Sigma \stackrel{q}{\to} G
\]
of an \'etale groupoid by the unit circle $\TT$.  The main difference here with Definition~\ref{DT} is that we assume $\TT$ has the usual topology, which is not discrete.  We can construct the line bundle for this extension as in Section~\ref{TSA}; the set of continuous compactly supported sections of the line bundle is a $*$-algebra under

\[
f*g(\gamma)=\sum_{r(\eta)=r(\gamma)}f(\eta)g(\eta\inv\gamma)\quad\quad{\text{and}}\quad\quad f^*(\gamma)=\overline{f(\gamma\inv)}.
\] 
The reduced $C^*$-algebra $C_r^*(G;\Sigma)$ is then the completion of the $*$-algebra of compactly supported sections with the reduced norm \cite{Ren08}.  Using \cite[Proposition~2.8]{BFPRny} (see also \cite[Proposition~II.4.2]{ren80}), there exists  an injective norm decreasing map from the reduced $C^*$-algebra of the twist to continuous sections of the line bundle which vanish at infinity that extends the identity map of compactly supported sections: $j: C_r^*(G;\Sigma)\to C_0(G;\Sigma)$.  We also note that $C_r^*(G;\Sigma)$ is nuclear if and only if $G$ is amenable \cite[Theorem~5.4]{Tak}.

\begin{thm}
\label{thm:C*}
Let $G$ be an amenable $\i2i$-groupoid Hausdorff locally compact groupoid, and let $\Sigma\to G$ be a twist.  Suppose $C$ is a C*-algebra such that 
\[
C_0(\go) \subseteq C \subseteq C^*_r(G;\Sigma).\]
Then $C_0(\go) \subseteq C$ is Cartan and hence $(C^*_r(G;\Sigma),C_0(\go))$ is purely Cartan.
\end{thm}

We will prove this theorem in a series of lemmas.  The work comes in showing that $(C,D)$ is a regular inclusion.
We first introduce a piece of notation.
We write the open support of $c \in C_r^*(G;\Sigma)$ as 
\[
S_c:= \{\gamma\in G : j(c(\gamma)) \neq 0\}.\]   



\begin{lem}
\label{lem:restrict}Let $G$ be a Hausdorff \'{e}tale groupoid and let $B\subseteq G$ be an open bisection and let $\gamma \in B.$  
Suppose $a \in C^*_r(G;\Sigma)$.  Then
there is a bisection $B' \subseteq B$ with $\gamma \in B'$, and a normalizer $n \in N(C_r^*(G;\Sigma),C_0(\go))$ such that $S_n \subseteq S_a$ and $j(n)(\sigma) = j(a)(\sigma)$ for all $\sigma \in B'$.

When $G$ is ample and $B$ is a compact open bisection, $n$ can be chosen so that $j(n) = j(a)|_B$.
\end{lem}

\begin{proof}
The key to the proof is that the reduced and uniform norms agree on functions supported on bisections, see for example \cite[Corollary 9.3.4]{SimsNotes}.
As $G$ is locally compact and Hausdorff, there is a nonempty open set $B' \subseteq B$ with $\gamma \in B'$, such that $B' \subseteq \overline{B'} \subseteq B$ and $\overline{B'}$ compact.
Let $h \in C_0(G)$ satisfy $h(\eta) = 1$ for $\eta \in \overline{B'}$, and $h(\eta) = 0$ for $\eta\in G \backslash B$.
Let $\{f_k\} \in C_c(G;\Sigma)$ be a sequence such that $f_k \to a$ in $\|\cdot\|_r$.  Then by the continuity of $j$ and the definition of the uniform norm we have
\[f_k=j(f_k) \to j(a), \]
in the uniform norm.
So $\{f_k\cdot h\}_k$ is Cauchy in the uniform norm, where $\cdot$ is being used to denote pointwise multiplication.
Further, each $f_k\cdot h \in C_c(G;\Sigma)$. Since the supports of each $f_k \cdot h$ are inside the bisection $B$, $\{f_k\cdot h\}_k$ 
is also Cauchy in the reduced norm.  Thus there exists $n \in C^*_r(G;\Sigma)$ such that 
$f_k\cdot h \to n$ in the reduced norm.
By the injectivity of $j$, we have that 
$j(n) = j(a)\cdot h$.
Hence $j(n)(\gamma) = j(a)(\gamma)$ for all $\gamma \in B'$.
The element $n$ is supported on the bisection $B \cap S_a$, and hence $n$ is a normalizer \cite[Proposition~4.8]{Ren08} with the desired properties.

When $G$ is ample and $B$ is a compact open bisection, we can let $B'=B$ and choose $h = 1_B$.
Then $n$ will satisfy $j(n) = j(a)|_B$.
\end{proof}

\begin{lem}
\label{lem: Cst normalizes} Let $G$ be an $\i2i$-groupoid that is locally compact, Hausdorff and amenable, and let $\Sigma\to G$ be a twist.  Suppose $a\in C^*_r(G;\Sigma)$ and $\gamma\in G$ with $j(a)(\gamma)\neq 0$.  Then there exists  $n\in N(C^*_r(G;\Sigma), C_0(\go))\cap C^*(\{a\}\cup C_0(\go))$ such that $j(n)(\gamma)\neq 0$.
\end{lem}

\begin{proof}
 Let $H:=(G \setminus \operatorname{Iso}(G)) \cup \go$.  Then $H$ is a principal wide subgroupoid of $G$.  Since $\go$ is open and $\operatorname{Iso}(G)$ is closed, $H$ is open.
 As $G$ is amenable and $H$ is a wide open subgroupoid.  Since any open set is locally closed, $H$ is also an amenable groupoid \cite[Proposition~5.1.1]{AnaRen00}.
Thus $C_0(\go) \subseteq C^*(H;\Sigma|_H)$ is Cartan by \cite[Theorem~3.3]{BEFPR}.
By \cite[Theorem~4.5]{BEFPR}, for every $C^*$-algebra $D$ such that $C_0(\go) \subseteq D \subseteq C^*(H; \Sigma|_H)$,  $C_0(\go) \subseteq D$ is Cartan.
Further, by \cite[Theorem~4.5]{BEFPR}, each such $D$ is of the form $D=C^*(L;\Sigma|_L)$ where $L$ is the wide open subgroupoid of $H$ defined by
\[L:= \{\gamma \in G \colon j(d(\gamma)) \neq 0 \text{ for some } d \in D\}.\]
Let  $C=C^*(\{a\}\cup C_0(\go))$ then
\[
C_0(\go) \subseteq C \subseteq C^*_r(G; \Sigma).\]
Moreover,
\[C_0(\go) \subseteq C \cap C^*(H; \Sigma|_H) \subseteq C^*(H; \Sigma|_H)\]
 is Cartan and hence 
$C \cap C^*(H;\Sigma|_H) = C^*(L;\Sigma|_L) \subseteq C$ where 
\[
L:= \{\gamma \in G \colon j(c(\gamma)) \neq 0 \text{ for some } c \in C \cap C^*(H; \Sigma|_H)\}.\]

By subtracting $E(a)$ from $a$, we can assume $S_a \subseteq G \setminus \go$. 
Fix $\gamma \in S_a$.  
Since $G$ is \'{e}tale and $S_a$ is open, there exists an open bisection $B$ such that $\gamma \in B \subseteq S_a$.
We consider two cases.  First, suppose $\gamma \in H$.
As $\gamma \in H$ we can choose $B$ so that $r(B)$ and $s(B)$ are disjoint.
Choose $f, g \in C_0(\go)$ so that $f$ is supported on $r(B)$, $g$ is supported on $s(B)$, and $f(r(\gamma))=g(s(\gamma)) = 1$.
For any $\sigma \in G$
\[ j(fag)(\sigma) = f(r(\sigma))j(a)(\sigma)g(s(\sigma)).\]
Since $r(B)$ and $s(B)$ are disjoint, the choice of $f$ and $g$ mean that 
\[ fag \in C \cap C_r^*(H; \Sigma|_H) = C_r^*(L; \Sigma|_L). \]
Further $j(fag)(\gamma) = j(a)(\gamma) \neq 0$, hence $\gamma \in L$.
If $B'$ is any bisection in $L$ with $\gamma \in B'$, and $n \in C_0(B'; \Sigma|_{B'})$ satisfying $n(\gamma) \neq 0$, then $n$ is a normalizer with the desired properties.

For the second case, suppose $s(\gamma)=r(\gamma) \neq \gamma$.
Let $B$ be an open bisection in $S_a$ with $\gamma \in B$.
And let $B' \subseteq B$ be a bisection with $\gamma \in B' \subseteq \overline{B'}\subseteq B$.
Choose $f, g \in C_0(\go)$ with $f$ supported on $r(B')$, $g$ supported on $s(B')$, and $f(r(\gamma))=g(s(\gamma))= 1$.
Write
\[\tilde{a}:=fag  \in C.\]
We claim that 
\[(G \setminus B') \cap S_{\tilde{a}} \subseteq L.\]
First, we show the left-hand side is in $H$.  Aiming for a contradiction, suppose not.  Then there exists $\alpha$ in the left-hand side with
$s(\alpha)=r(\alpha) \neq \alpha$.  Since $G$ is $\i2i$, $Gs(\alpha)=r(\alpha) G=r(\alpha)Gs(\alpha)=\{\alpha, s(\alpha)\}$.  Since $\alpha\in (G \setminus B') \cap S_{\tilde{a}}  \subseteq S_{\tilde{a}}$, we have  $s(\alpha) \in s(B')$. So there exists $\eta \in B'\subseteq S_a\subseteq G-\go$ with $s(\eta)=s(\alpha)$.  But then $\eta\in Gs(\alpha)\backslash \go= \{\alpha\}$, that is, $\eta=\alpha$ which is a contradiction,
proving the left-hand-side is in $H$. 
Arguing as in case 1 then shows that $(G\setminus B')\cap S_{\tilde{a}} \subseteq L$.

By Lemma~\ref{lem:restrict}, we can choose $n \in N(C_r^*(G;\Sigma),C_0(\go))$ such that $S_n \subseteq S_{\tilde{a}}$, and $j(n)(\sigma) = j(\tilde{a})(\sigma)$ for all $\sigma \in B'$.
Let $ e= \tilde{a} - n$.
Notice 
\[S_e \subseteq (G \setminus B') \cap S_{\tilde{a}} \subseteq L\]
and so $e \in C^*(L;\Sigma|_L) \subseteq C$.  But then $n \in C$ as well, and is a normalizer as needed. 
\end{proof}

\begin{proof}[Proof of Theorem~\ref{thm:C*}]
The proof follows almost exactly as the proof of \cite[Theorem~4.5]{BEFPR} with Lemma~\ref{lem: Cst normalizes} used in place of \cite[Lemma~4.3]{BEFPR}.\end{proof}

\begin{rmk}
The counter examples in Section~\ref{PqC} translate to the $C^*$-setting virtually unchanged.  For Proposition~\ref{prop: two arrows}, pick sections $\phi_i$ with support in $\Ooo_i$ and take $f=\phi_1+\phi_2$  and $C=\overline{\spn}\{g_1 f g_2: g_i\in C_0(\go)\}$ in the proof.  For  Lemma~\ref{Bad apple} and Proposition~\ref{bad apple}, it is only necessary to change the twisted Steinberg algebras to twisted groupoid $C^*$-algebras.  

\end{rmk}

In light of the above remark we get the following corollary.

\begin{corollary}\label{cor pC cstar} Let $G$ be an amenable Hausdorff \'etale groupoid.  Then $(C_r^*(G), C_0(\go))$ is purely Cartan if and only if $G$ is $\i2i$.\end{corollary}

\begin{rmk}\label{rmk: KenUrs}
    In a different, but related, setting Kennedy and Ursu \cite{KenUrs} have recently proved a result similar to Theorem~\ref{thm:C*}.
    Let $\Gamma$ be a discrete group acting on a C$^*$-algebra $B$.
    They study when intermediate algebras $B \subseteq C \subseteq B\rtimes_r \Gamma$ come from subactions of $\Gamma$ on $B$.
    In the case when $B = C(X)$ is abelian and $\Gamma$ is amenable, their results are a special case of what is shown here.
\end{rmk}

\appendix
\section{Conditional expectations } \label{CE}

We establish another characterization of the special conditional expectations $\ex: A_R(G,\Sigma)\to A_R(\go)$ when $G$ is principal and $\Char{R}\neq 2.$
This is inspired by the $C^*$-algebra case as described in the following remark.

\begin{rmk}
Let $A$ be a unital $C^*$-algebra and let $D \subseteq A$ be a maximal abelian C$^*$-subalgebra. Archbold, Bunce and Gregson \cite[Corollary~2.7]{ABG} showed the following equivalent:
\begin{enumerate}
    \item every pure state on $D$ extends uniquely to a pure state on $A$;
    \item\label{c*ii} there is a faithful conditional expectation $\ex\colon A \rightarrow D$ satisfying \[\{\ex(a)\} = \overline{\mathrm{conv}}\{uau^* \colon u \in D\text{ unitary}\} \cap D.\]
\end{enumerate}
A special case of this result was proved by Power \cite{Pow} for crossed products $C(X) \rtimes_r \ZZ$ by free actions.
Phillips later generalized (and simplified) Power's approach to study free actions by discrete groups \cite{GKPT}. 
\end{rmk}

We establish two technical lemmas before proving Theorem~\ref{thm: 1} about conditional expectations in the principal setting.
The first lemma is a refinement of Proposition~\ref{prop regular} for Steinberg algebras of principal groupoids.
In the proof of Lemma~\ref{lem: 1}, we adapt Phillips' approach to groupoids and Steinberg algebras, see \cite[Proposition~11.1.19]{GKPT}.

\begin{lem}\label{lem: 2}
Suppose $R$ is a commutative ring with identity and $(\Sigma,i,q)$ is a discrete twist  over an ample Hausdorff groupoid $G$, with $A_R(\go)$  $\WT$ with respect to $R$.
Suppose that $G$ is principal. 
If $f \in \Stein$, then $f$ can be written as  $f = \sum_{i=1}^N r_i \chrp{B_i}$ where
\begin{enumerate}
    \item the sets $B_i \subseteq G$ together form a disjoint collection of compact open bisections;
    \item either $B_i \subseteq G^{(0)}$ or $B_i \cap G^{(0)} = \emptyset$;
    \item if $B_i \cap G^{(0)} = \emptyset$ then $r(B_i) \cap s(B_i) =  \emptyset$ for each $1\leq i \leq N$.
\end{enumerate}
\end{lem}

\begin{proof}
By Proposition~\ref{prop regular} (\cite[Proposition~2.8]{A+8}), we can write $f$ as $f = \sum_{i=1}^N r_i \chrp{B_i}$ where $\{B_1,\ldots, B_n\}$ are mutually disjoint compact open bisections in $G$ and $\psi_{B_i}: R\times B_i\to \lineP\inv(B_i)$ are $R$-linear homeomorphisms.
Hence condition (1) holds. 
To show that conditions (2) and (3) hold it suffices to prove the result for $f = \chrp{B}$ where $B$ is a compact open bisection in $G$ and $\psi_{B}: R\times B \to \lineP\inv(B_i)$ is a  $R$-linear homeomorphism.
Note that $\go$ is both open and closed since $G$ is \'{e}tale and Hausdorff. 
Thus, if $B\subseteq G$ is a compact open bisection it can be split into two disjoint compact open bisections: $B \cap \go$ and $B - \go$ and hence (2) holds.
To show (3), we can assume $B \subseteq G \backslash G^{(0)}$. Since $G$ is principal, each $\gamma \in B$ has $r(\gamma) \neq s(\gamma)$. Thus, using that $\go$ is Hausdorff, we can find a compact open bisection $B_\gamma \subseteq B$ containing $\gamma$ satisfying $s(B_\gamma) \cap r(B_\gamma) = \emptyset.$  As $B$ is compact, we can take a finite subset of $\{B_\gamma\}_{\gamma \in B}$ covering $B$. Using the disjointification of this finite subcover \cite[Remark~2.4]{CFST}, write $B$ as
\[B = \bigcup_{i=1}^n B_i \]
where
\begin{enumerate}
    \item the sets $B_i$ together form a disjoint collection of compact open bisections that cover of $B$;
    \item $r(B_i) \cap s(B_i) =  \emptyset$ for each $1\leq i \leq n$.
\end{enumerate}
The result follows.
\end{proof}

\begin{lem}\label{lem: 1}
Suppose $R$ is a commutative ring with identity and $(\Sigma,i,q)$ is a discrete twist  over a  Hausdorff ample groupoid $G$ with $A_R(\go)$  $\WT$ over $R$.
Suppose $B_1, \ldots, B_k$ are bisections such that $r(B_i) \cap s(B_i) = \emptyset$.
Then for any compact $K \subseteq G^{(0)}$ such that $\bigcup_{i=1}^k s(B_i) \cup \bigcup_{i=1}^k r(B_i) \subseteq K$, there are $u_1,\ldots, u_{2^k} \in A_R(G^{(0)})$ such that
\begin{enumerate}
\item $\sum_{j=1}^{2^k} u_j(r(\beta))u_j(s(\beta)) = 0,$ for all $\beta \in \bigcup_{i=1}^k B_i$;
\item $\IM(u_i|_K)\in \{-1,1\}$, for $1\leq i \leq 2^k$.
\end{enumerate}
\end{lem}

\begin{proof}
Take a compact open bisection $B$ which satisfies $r(B) \cap s(B) = \emptyset,$ and let $K \subseteq G^{(0)}$ be a compact set so that $s(B)\cup r(B) \subseteq K$.
As $G^{(0)}$ is locally compact, Hausdorff and totally disconnected there is an open and compact $L\subseteq G^{(0)}$ such that $K \subseteq L$.
Let $u_1 = \chrp{L}$ and let $u_2 = \chrp{L}- 2 \chrp{r(B)}$.
Then $u_1,u_2$ satisfy  $\IM(u_1|_K)\cup\IM(u_2|_K) \in \{-1,1\}$.  Also, for any $\beta \in B$,

\[
u_1(r(\beta))u_1(s(\beta)) + u_2(r(\beta))u_2(s(\beta))  = (1)(1) + (-1)(1) = 0. 
\]
Thus, the result holds when $k=1$.

Now fix  bisections $B_1,\ldots, B_k$ such that  $r(B_i) \cap s(B_i) = \emptyset$ for each $1\leq i \leq k$; and a compact set $K$ such that $\bigcup_{i=1}^k s(B_i) \cup \bigcup_{i=1}^k r(B_i) \subseteq K$.
Assume the result holds for $l< k$.
By hypothesis, there are $u_1, \ldots, u_{2^{k-1}} \in A_R(G^{(0)})$ with $\IM(u_i|_K)\in \{-1,1\}$,  for $1\leq i \leq 2^{k-1}$ and so that 
\[
 \sum_{j=1}^{2^{k-1}} u_j(r(\beta))u_j(s(\beta)) = 0,
 \]
for all $\beta \in \bigcup_{i=1}^{k-1} B_i$.
Also, by the above argument for $k=1$, there are $w_1, w_2 \in A_R(G^{(0)})$ with $\IM(w_1|_K)\cup\IM(w_2|_K) \in \{-1,1\}$ so that 
\[
w_1(r(\beta))w_1(s(\beta)) + w_2(r(\beta))w_2(s(\beta))= 0,
\]
for any $\beta \in B_k$.

Now take $\beta \in \bigcup_{i=1}^k B_i$. We have that
\begin{align*}
    \sum_{l=1}^2\sum_{j=1}^{2^{k-1}} (w_l u_j)(r(\beta))(w_l u_j)(s(\beta)) &=
    \left(\sum_{l=1}^2 w_l(r(\beta))w_l(r(\beta))\right)
    \left(\sum_{j=1}^{2^{k-1}} u_j(r(\beta))u_j(s(\beta))\right)\\& = 0.
\end{align*}
Thus the result holds using the $2^k$ elements  $w_l u_j\in A_R(G^{(0)})$, $1\leq l \leq 2$, $1\leq j \leq 2^{k-1}$.
\end{proof}

We are now ready to give an algebraic description of the conditional expectation $\ex \colon \Stein \rightarrow A_R(\go)$ when $R$ is an integral domain and $G$ is principal.

\begin{thm}\label{thm: 1}
Let $R$ be an integral domain, and let $(\Sigma,i,q)$ be a discrete twist  over a  Hausdorff ample groupoid $G$ with $A_R(\go)$  $\WT$ with respect to $R$.
Let $\ex(f) = f|_{G^{(0)}}$ be the conditional expectation from  $\Stein$ onto $A_R(G^{(0)})$.

\begin{enumerate}
\item Assume $\Char(R)\neq 2$. If $G$ is principal, then for each $f \in A_R(G)$ there exist $u_1, \ldots, u_{2^k} \in A_R(G^{(0)})$ such that
$$\ex(f) = \frac{1}{2^k} \sum_{i=1}^{2^k}  u_i * f * u_i.$$

\item If for each $f \in \Stein$ there exist $u_1, \ldots, u_N \in A_R(G^{(0)})$ such that
$$\ex(f) = \frac{1}{N} \sum_{i=1}^N u_i * f * u_i,$$
then $G$ is principal.
\end{enumerate}
\end{thm}

\begin{proof}

Take $f \in \Stein$.
If the support $S_f \subseteq G^{(0)}$ then we can take $n=1$ and $u_1 = 1_{S_f}$ and we are done.
Thus we will assume that $S_f \not\subseteq G^{(0)}$.
By Lemma~\ref{lem: 2}, we can write $f$ as

\[
f = \sum_{i=1}^m \re_i \chrp{C_i} + \sum_{i=1}^k \rea_i \chrp{B_i}, 
\]
where $\{C_i\}_i$ is a collection of disjoint compact bisections in $G^{(0)}$ and $\{B_i\}_i$ is a collection of disjoint compact bisections in $G\backslash G^{(0)}$ satisfying $r(B_i) \cap s(B_i) = \emptyset$ for $1\leq i \leq k$,  $\psi_{B_i}: R\times B_i\to \lineP\inv(B_i)$ are $R$-linear homeomorphisms and $\re_i,\rea_i\in R$.

Thus $S_f = \bigcup_{i=1}^m C_i \cup \bigcup_{i=1}^k B_i$. Let
\[
 K = \bigcup_{i=1}^m C_i \cup \bigcup_{i=1}^k(r(B_i) \cup s(B_i)). 
 \]
By Lemma~\ref{lem: 1} there are $u_1,\ldots,u_{2^k} \in A_R(G^{(0)})$ so that $\IM(u_i|_K)\in \{-1,1\}$ and
$$ \sum_{j=1}^{2^k} u_j(r(\beta))u_j(s(\beta)) = 0,$$
for all $\beta \in \bigcup_{i=1}^k B_i$.
We now have, for any $\gamma\in G$,
\begin{align*}
    \sum_{i=1}^{2^k} u_i * f * u_i(\gamma) &= \sum_{i=1}^{2^k} u_i(r(\gamma)) f(\gamma) u_i(s(\gamma))\\
    &= \left(\sum_{i=1}^{2^k} u_i(r(\gamma))u_i(s(\gamma))\right)f(\gamma).
\end{align*}
Since $S_f =  \bigcup_{i=1}^m C_i \cup \bigcup_{i=1}^k B_i$, we thus have
\begin{align*}
   \sum_{i=1}^{2^k} u_i * f * u_i(\gamma) &= \begin{cases}0  &\text{if } \gamma \notin S_f, \\
    0 & \text{if } \gamma \in S_f \backslash G^{(0)}, \\
    2^k f(\gamma) &\text{if }\gamma \in S_f \cap G^{(0)}.\end{cases}
\end{align*}
Thus
$$\sum_{i=1}^{2^k} u_i * f * u_i = 2^k \sum_{i=1}^m \re_i 1_{C_i} = 2^k \ex(f). $$

Now if $2^k \ex(f) =0$, then $2^k$ is a zero divisor in $R$ or $\ex(f) = 0.$
If $2^k$ is a zero divisor in $R$, then there exists a nonzero $r \in R$ such that $2^k r = 0$.
It follows that the order of $r$ divides $2^k$.
By assumption, $\ord(r)$ does not divide $2$, and hence $2^k\ex(f) = 0$ if and only if $\ex(f) = 0$.

Recall that $A_R(\go)$ is  $\WT$ over $R$ and $A_r(\go)$ is generated by its idempotents.
This implies $A_R(\go)$ is torsion-free over $R$, see \cite[Section~2]{A+8}.
As we are assuming that $\Char(R)$ is not divisible by $2$, we have that if $g\in A_R(G^{(0)})$, then $2^k g = 0$ if and only if $g = 0.$
It follows that if $h,g \in A_R(G^{0})$, then $2^kh = 2^kg$ if and only if $h = g.$
Thus
\[
 \ex(f) = \frac{1}{2^k}\sum_{i=1}^{2^k} u_i * f * u_i. 
 \]

Suppose now that $G$ is not principal. 
Choose $\gamma \in G\backslash G^{(0)}$ such that $r(\gamma) = s(\gamma)$.
Take any $f \in \Stein$ such that $f(\gamma)\neq 0$ and $f(r(\gamma))\neq 0$. 
For example, for $B$ a compact open bisection containing $B$ we can take $f=\chrp{B}+\chrp{r(B)}$, where $\psi_{B}: R\times B\to \lineP\inv(B_i)$  and $\psi_{r(B)}: R\times r(B)\to \lineP\inv(r(B))$ are $R$-linear homeomorphisms.

Then, for any  $u_1,\ldots,u_N\in A_R(G^{(0)})$ we have 
\begin{align*}
    \sum_{i=1}^N u_i * f * u_i(\gamma) &= \sum_{i=1}^n u_i(r(\gamma)) f(\gamma) u_i(s(\gamma)) \\  
    &=\left(\sum_{i=1}^N |u_i(r(\gamma))|^2\right) f(\gamma),
\end{align*}
and 
\begin{align*}
    \sum_{i=1}^N u_i * f * u_i(r(\gamma)) &= \sum_{i=1}^n u_i(r(\gamma)) f(r(\gamma)) u_i(r(\gamma)) \\  
    &=\left(\sum_{i=1}^N |u_i(r(\gamma))|^2\right) f(r(\gamma)).
\end{align*}
Hence $\sum_{i=1}^N u_i * f * u_i(\gamma) = 0$ if and only if $\sum_{i=1}^N u_i * f * u_i(r(\gamma)) = 0$.
Thus $\ex(f) \neq \frac{1}{N}\sum_{i=1}^n u_i * f * u_i$.
\end{proof}


\begin{rmk}
We can prove Lemma~\ref{normalizers in C} by finding a normalizer in $\bi(c)$ for a given element $c\in \Steink$ using Theorem~\ref{thm: 1} in the case $\Char(K)\neq 2$ as follows.
Recall that Lemma~\ref{normalizers in C} is the key step in the proofs of both Theorem~\ref{thm: diagonal case} and Theorem~\ref{STbi}.

By Lemma~\ref{lem: 1}, we can write
$$ c = \sum \re_i \chrp{B_i}$$
where $B_1,\ldots,B_m$ are disjoint open compact bisections satisfying either $B_i \subseteq \go$ or $r(B_i) \cap s(B_i) = \emptyset$ for $i = 1,\ldots, m$, $\psi_{B_i}: R\times B_i\to \lineP\inv(B_i)$ are $R$-linear homomorphisms, and $\re_i\in R$. There is a unique $i_0$ such that $\gamma \in B_{i_0}$.
Let $n = \chrp{B_{i_0}}$.
Then $c|_{S_n} = \re_{i_0}n$.

Let $\ex$ denote the conditional expectation from $\Steink$ to $A_K(\go)$.
By Theorem~\ref{thm: 1} there are $u_1,\ldots,u_{2^k} \in A_K(\go)$ so that
\begin{align*}
    n \ex(n^\dag c) &= n \left(\frac{1}{2^k}\sum_{j=1}^{2^k}u_j n^\dag c u_j\right) \\
    &= \sum_{j=1}^{2^k} \left(\frac{n u_j n^\dag}{2^k}\right) c u_j.
\end{align*}
As $n$ normalizes $A_K(\go)$ it follows that $n \ex(n^\dag c) \in \bi(c)$.
Using the fact that the bisections $B_1,\ldots, B_m$ are disjoint, a straightforward calculation shows that $n \ex(n^\dag c) = \alpha_{i_0} n$, and Lemma~\ref{normalizers in C} follows.
\end{rmk}

\end{document}